\documentclass[letterpaper,oneside,10pt]{article}
\usepackage[USenglish]{babel} %francais, polish, spanish, ...
\usepackage[ansinew]{inputenc}
\usepackage{verbatim}
\usepackage{amsmath, amsthm, amssymb}

\usepackage{graphicx,color}
\usepackage{float}
\usepackage{latexsym}
\usepackage{amsmath}
\usepackage{amsfonts}
\usepackage{enumerate}
\usepackage{fancyhdr}
\usepackage{fullpage}
\renewcommand{\Re}{\ensuremath{\mathbb{R}}}
\begin{document}
\newcommand{\sg}{\mathbf{s}_g}
\newcommand{\bg}{\mathbf{b}_g}
\newcommand{\sa}{\mathbf{s}_a}
\newcommand{\ba}{\mathbf{b}_a}
\newcommand{\sgh}{\hat{\mathbf{s}}_g}
\newcommand{\bgh}{\hat{\mathbf{b}}_g}
\newcommand{\sah}{\hat{\mathbf{s}}_a}
\newcommand{\bah}{\hat{\mathbf{b}}_a}
\newcommand{\mw}{\left\|\omega\right\|}
\newcommand{\dd}{\partial}
%%\pagestyle{empty} %No headings for the first pages.

%%%%%SET HEADER AND FOOTER%%%%%%%%%%%%%%%%%%%%%%%%%%%%%

%\fancyhead{}
%\fancyfoot{}
%\fancyhead[RO,RE]{DACES Status Report 2\\AFOSR-STTR FA9550-10-C-0039\\AR10-107\\October 30, 2010}
%\fancyfoot[CO,CE]{AURORA PROPRIETARY}
%\fancyfoot[RO, LE] {Page \thepage}
%\fancyhead[L]{ % right
%   \includegraphics[height=0.53in]{./figs/Aurora_logo.png}
%}

%\newpage
\pagestyle{fancyplain}
\renewcommand{\headrulewidth}{0.2pt}
\renewcommand{\footrulewidth}{0.2pt}

%\newpage
%
%%% Title Page %%%%%%%%%%%%%%%%%%%%%%%%%%%%%%%%%%%%%%%%%%%%%%%
%%\pagestyle{plain}
%\title{Model Reference Adaptive Control with Input and Output Constraints\\}
%\author{\\ Mehrdad Pakmehr \\ Georgia Institute of Technology \\ School of Aerospace Engineering}
%
%\maketitle
%\newpage
%%\tableofcontents%Table of contents
%%\pagestyle{fancyplain} %Now display headings: headings / fancy / ...
%\footskip = 15pt
%
%\newpage

%% The class has several options
%  onecolumn/twocolumn - format for one or two columns per page
%  10pt/11pt/12pt - use 10, 11, or 12 point font
%  oneside/twoside - format for oneside/twosided printing
%  final/draft - format for final/draft copy
%  cleanfoot - take out copyright info in footer leave page number
%  cleanhead - take out the conference banner on the title page
%  titlepage/notitlepage - put in titlepage or leave out titlepage
%
%% The default is oneside, onecolumn, 10pt, final

\pagestyle{plain}
\title{\textbf{Gain Scheduling Control of Gas Turbine Engines: Absolute Stability by Finding a Common Lyapunov Matrix}\\}

\author{\\ Mehrdad Pakmehr \footnote{PhD Candidate, School of Aerospace Engineering, Georgia Institute of Technology, email: mehrdad.pakmehr@gatech.edu.},
Nathan Fitzgerald \footnote{Propulsion Development Engineer, Aurora Flight Sciences Corporation.},
Eric Feron \footnote{Professor, School of Aerospace Engineering, Georgia Institute of Technology.},
Jeff S. Shamma \footnote{Professor, School of Electrical and Computer Engineering, Georgia Institute of Technology.},
Alireza Behbahani \footnote{Senior Aerospace Engineer, Air Force Research Laboratory (AFRL).}\\}

%\author{Mehrdad Pakmehr\thanks{Address all correspondence to this author.}
%    \affiliation{
%		PhD Candidate\\
%		Decision and Control Laboratory (DCL)\\
%		School of Aerospace Engineering\\
%		Georgia Institute of Technology\\
%		Atlanta, Georgia 30332\\
%  	Email: mehrdad.pakmehr@gatech.edu   }	 }
%
%\author{Nathan Fitzgerald
%    \affiliation{ Propulsion Development Engineer\\
%		Aurora Flight Sciences Corporation\\
%		Manassas, VA 20110\\
%    Email: nfitzgerald@aurora.aero  } }
%
%\author{Eric M. Feron
%     \affiliation{Professor\\
%    School of Aerospace Engineering\\
%	Georgia Institute of Technology\\
%		Atlanta, Georgia 30332\\
%    Email: feron@gatech.edu  } }
%
%\author{Jeff S. Shamma
%    \affiliation{Professor\\
%    School of Electrical and Computer \\ Engineering\\
%	Georgia Institute of Technology\\
%		Atlanta, Georgia 30332\\
%    Email: shamma@gatech.edu  } }
%
%\author{Alireza Behbahani
%    \affiliation{Senior Aerospace Engineer\\
%    Air Force Research Laboratory\\
%    Wright-Patterson Air Force Base, Ohio 45433\\
%    Email: alireza.behbahani@wpafb.af.mil }  }

\date{\null}

\maketitle
\footskip = 15pt

%%%%%%%%%%%%%%%%%%%%%%%%%%%%%%%%%%%%%%%%%%%%%%%%%%%%%%%%%%%%%%%%%%%%%%
%\begin{multicols}{2}

\section*{Abstract}

This manuscript aims to develop and describe gain scheduling control concept for a gas turbine engine which drives a variable pitch propeller. An architecture for gain-scheduling  control is developed that controls the turboshaft engine for large thrust commands in stable fashion with good performance. Fuel flow and propeller pitch angle are the two control inputs of the system. New stability proof has been developed for gain scheduling control of gas turbine engines using global linearization and LMI techniques. This approach guarantees absolute stability of the closed loop gas turbine engines with gain-scheduling controllers.

%%%%%%%%%%%%%%%%%%%%%%%%%%%%%%%%%%%%%%%%%%%%%%%%%%%%%%%%%%%%%%%%%%%%%%
\section*{Nomenclature}

$N_1$: Non-dimensional Fan Spool Speed\\
$N_2$: Non-dimensional Core Spool Speed\\
$T$: hrust (N)\\
TSFC: Thrust Specific Fuel Consumption\\
$\alpha$: Scheduling Parameter\\
$\sigma$: Singular Value\\
$\lambda$: Eigenvalue\\

%%%%%%%%%%%%%%%%%%%%%%%%%%%%%%%%%%%%%%%%%%%%%%%%%%%%%%%%%%%%%%%%%%%%%%
\section{Introduction}

The gain-scheduling approach is perhaps one of the most popular nonlinear control design approaches which has
been widely and successfully applied in fields ranging from aerospace to process control \cite{research-rugh-2000, surveyGS-leith-2000}. Gain-scheduling, specifically has been used for gas turbine engine control, some of these works are \cite{lpv-balas-2002, lpv-gilbert-2010, lpv-bruzelius-2002, lpv-Shuqing-2010, approximate-zhao-2011, approximate-yu-2011}. In general, stability and control of gas turbine engines have been of interest to researchers and engineers from a variety of perspectives. Stability of axial flow fans operating in parallel has been investigated in \cite{fanStability-simon-1985}. An application of robust stability analysis tools for uncertain turbine engine systems is presented in \cite{robustAeroengine-arriffin-1997}. Application of the Linear-Quadratic-Gaussian with Loop-Transfer-Recovery methodology to design of a control system for a simplified turbofan engine model is considered in \cite{lqg-garg-1989}. A unified robust multivariable approach to propulsion control design has been developed in \cite{turbofanControl-fredrick-2000}. A simplified scheme for scheduling multivariable controllers for robust performance over a wide range of turbofan engine operating points is presented in \cite{turbofanSched-garg-1997}.

In the previous work by authors \cite{distributed-pakmehr-2009, decentralized-pakmehr-2010} controllers developed for single spool and twin spool turboshaft system. Those controllers were designed for small transients, and small throttle commands. In this work we develop a gain-scheduling control structure for JetCat SPT5 turboshaft engine using the method presented in \cite{gainsched-shamma-1988, research-rugh-2000, gain-shamma-2006, overview-shamma-2012}. the controller is designed to be used for entire flight envelope of the twin spool turboshaft engine.

In this manuscript, first a linear representation of the turbofan system dynamics is developed. Then control theoretic concepts for gain-scheduling control of this model is presented. The developed controller can be used for the entire flight envelope of the engine with guaranteed stability. Finally the simulation results for gain scheduling control of a physics-based nonlinear model of the JetCat SPT5 turboshaft engine are presented.

%%%%%%%%%%%%%%%%%%%%%%%%%%%%%%%%%%%%%%%%%%%%%%%%%%%%%%%%%%%%%%%%%%%%%%
\section{Gain Scheduling Control Design}

Consider the nonlinear dynamical system
\begin{equation}\label{eqn_gs1}
\begin{array}{c}
\dot{x}^p(t)= f^p(x^p(t),u(t)),\\[5pt]
y(t)=g^p(x^p(t),u(t)),
\end{array}
\end{equation}
where $x^p \in \Re^n$ is the state vector, $u\in \Re^m$ is the control input vector, $y\in \Re^k$ is the output vector, $f^p(.)$ is an $n$-dimensional, and $g^p(.)$ is an $k$-dimensional differentiable nonlinear vector functions.
We want to design a feedback control such that $y(t) \rightarrow r(t)$ as  $t \rightarrow \infty$, where $r(t) \in D_r \subset \Re^k$ is the output reference signals vector.

Assume that for each $r(t) \in D_r$, there is a unique pair $(x^p_e, u_e)$ that depends continuously on $r$ and satisfies the equations:

\begin{equation}\label{eqn_gs2}
\begin{array}{c}
0= f^p(x^p_e,u_e),\\[5pt]
r=g^p(x^p_e,u_e),
\end{array}
\end{equation}
in case of a constant $r$. $x^p_e$ is the desired equilibrium point and $u_e$ is the steady-state control that is needed to maintain equilibrium at $x^p_e$.

\newtheorem{deff}{Definition}
\begin{deff} \label{def1}
The functions $x^p_e(\alpha), u_e(\alpha)$, and $r_e(\alpha)$ define an equilibrium family for the plant (\ref{eqn_gs1}) on the set $\Omega$ if
\begin{equation}\label{eqn_gs222}
\begin{array}{l}
f^p(x^p_e(\alpha),u_e(\alpha), r_e(\alpha))=0, \\[5pt]
g^p(x^p_e(\alpha),u_e(\alpha))=r_e(\alpha),  ~\alpha \in \Omega.
\end{array}
\end{equation}
\end{deff}

Let $\Omega \subset \Re^{m+n}$ be the region of interest for all possible system state and control vector $(x^p,u)$ during the system operation, and denote $x^{p*}_i$ and $u^*_i$, $i\in I = {1, 2, . . . , l}$, as a set of (constant) operating points located at some representative (and properly separated) points inside $\Omega$. Introduce a set of $l$ regions $\Omega_i$ centered at the chosen operating points $(x^{p*}_i, u^*_i)$, and denote their interiors as $\Omega_{i0}$, such that $\Omega_{j0} \bigcap \Omega_{k0}={\oslash}$ for all $j \neq k$, and $\bigcup_{i=1}^{l} \Omega_i=\Omega$. The linear model around each equilibrium point is
\begin{equation}\label{eqn_gs3}
\begin{array}{l}
\dot{x}^p = A^p_i (x^p-x_i^{p*}) + B^p_i (u-u_i^*),\\
y = C^p_i (x^-x_i^{p*}) + D^p_i (u-u_i^*) + y_i^* ,
\end{array}
\end{equation}
where the matrices are obtained as follows
\begin{equation}\label{eqn_gs4}
\begin{array}{l}
\displaystyle A^p_i=\frac{\partial f^p}{\partial x^p}|_{(x^{p*}_i, u^*_i)}, ~~ \forall (x^p,u) \in \Omega_i, \\[5pt]
\displaystyle B^p_i=\frac{\partial f^p}{\partial u}|_{(x^{p*}_i, u^*_i)}, ~~ \forall (x^p,u) \in \Omega_i, \\[5pt]
\displaystyle C^p_i=\frac{\partial g^p}{\partial x^p}|_{(x^{p*}_i, u^*_i)}, ~~ \forall (x^p,u) \in \Omega_i, \\[5pt]
\displaystyle D^p_i=\frac{\partial g^p}{\partial u}|_{(x^{p*}_i, u^*_i)}, ~~ \forall (x^p,u) \in \Omega_i.
\end{array}
\end{equation}
Here we assume that the common boundary of two regions $\Omega_j$ and $\Omega_z$ belongs to only one of $\Omega_j$ and $\Omega_z$.  Note that at each moment, $(x^p,u)$ belongs to only one $\Omega_i$.

Performing linearizations at a series of trim points gives a linearization family described by
\begin{equation}\label{eqn_gs5}
\begin{array}{l}
\delta \dot{x}^p = A^p(\alpha) \delta x^p + B^p(\alpha) \delta u,\\[5pt]
\delta y = C^p(\alpha) \delta x^p + D^p(\alpha) \delta u.
\end{array}
\end{equation}

where
\begin{equation}\label{eqn_gs6}
\begin{array}{l}
\delta x^p = x^p-x^p_e(\alpha) \\[5pt]
\delta y = y-y_e(\alpha),\\[5pt]
\delta u = u-u_e(\alpha),  ~~~ \forall \alpha \in \Omega.
\end{array}
\end{equation}

Gain scheduled controller for plant (\ref{eqn_gs5}), is designed as follows. First, a set of parameter values $\alpha_i$ are selected, which represent the range of the plant's dynamics, and a linear time-invariant controller for each is designed. Then, in between operating points, the controller gains are linearly interpolated such that for all frozen values of the parameters, the closed loop system has excellent properties, such as nominal stability and robust performance. To guarantee that the closed loop system will retain the feedback properties of the frozen-time designs, the scheduling variables should vary slowly withe respect to the system dynamics \cite{gainsched-shamma-1988}.

The parameter $\alpha$ called the scheduling variable in gain scheduling control. $A^p(\alpha), B^p(\alpha), C^p(\alpha)$, and $D^p(\alpha)$ are the parameterized linearization system matrices and $x^p_e(\alpha), u_e(\alpha)$, and $y_e(\alpha)$ are the parameterized steady-state system variables, which form the equilibrium manifold of system (\ref{eqn_gs1}). The subscript $'e'$ stands for steady-state throughout this paper.

Let the controller have the following structure
\begin{equation}\label{eqn_gs7}
\begin{array}{l}
\dot{x}^c = A^c(\alpha) \delta x^c + B^c(\alpha) [\delta y-\delta r],\\[5pt]
\delta u = C^c(\alpha) \delta x^c + D^c(\alpha) [\delta y-\delta r], ~~~ \forall \alpha \in \Omega.
\end{array}
\end{equation}

where
\begin{equation}\label{eqn_gs8}
\begin{array}{l}
\delta x^c = x^c-x^c_e(\alpha) \\[5pt]
\delta r = r-r_e(\alpha), ~~~ \forall \alpha \in \Omega.
\end{array}
\end{equation}

A standard realization of the parameterized controller can be written in the following form with the reference signal explicitly displayed

\begin{equation} \label{eqn_gs9}
\begin{array}{l}
  \left[
       \begin{array}{c}
           \dot{x}^c \\
           \delta u
       \end{array}
    \right] =
     \left[
       \begin{array}{ccc}
           A^c(\alpha) & B^c(\alpha) & -B^c(\alpha) \\
           C^c(\alpha) & D^c(\alpha) & -D^c(\alpha)
       \end{array}
    \right] ~
     \left[
       \begin{array}{c}
           \delta x^c \\
           \delta y \\
           \delta r
       \end{array}
    \right], ~~~ \forall \alpha \in \Omega.
\end{array}
\end{equation}

Controller has the general form
\begin{equation}\label{eqn_gs10}
\begin{array}{c}
\dot{x}^c(t)= f^c(x^c(t),y(t), r(t)),\\[5pt]
u(t)=g^c(x^c(t),y(t), r(t)),
\end{array}
\end{equation}
with the input and output signals corresponding to the nonlinear plant (\ref{eqn_gs1}).

The objective in linearization scheduling is that the equilibrium family of the controller (\ref{eqn_gs10}) match the plant equilibrium family, so that the closed loop system maintains suitable trim values, and second the linearization family of the controller is the designed family of linear controllers \cite{research-rugh-2000}.

For the equilibrium conditions of plant (\ref{eqn_gs1}) and controller (\ref{eqn_gs10}) to match, there must exist a function $x^c_e(\alpha)$ such that
\begin{equation}\label{eqn_gs11}
\begin{array}{l}
0= f^c(x^c_e(\alpha),y_e(\alpha),r_e(\alpha)),\\[5pt]
u_e(\alpha)=g^c(x^c_e(\alpha),y_e(\alpha), r_e(\alpha)), ~~~ \forall \alpha \in \Omega,
\end{array}
\end{equation}
where
\begin{equation}\label{eqn_gs12}
\begin{array}{l}
\displaystyle A^c(\alpha)=\frac{\partial f^c}{\partial x^c}|_{(x^c_e(\alpha),y_e(\alpha),r_e(\alpha))}, \\[5pt]
\displaystyle B^c(\alpha)=\frac{\partial f^c}{\partial y}|_{(x^c_e(\alpha),y_e(\alpha),r_e(\alpha))}, \\[5pt]
\displaystyle C^c(\alpha)=\frac{\partial g^c}{\partial x^c}|_{(x^c_e(\alpha),y_e(\alpha),r_e(\alpha))}, \\[5pt]
\displaystyle D^c(\alpha)=\frac{\partial g^c}{\partial y}|_{(x^c_e(\alpha),y_e(\alpha),r_e(\alpha))}, ~~~ \forall \alpha \in \Omega.
\end{array}
\end{equation}

So the controller family has the form
\begin{equation}\label{eqn_gs13}
\begin{array}{l}
\dot{x}^c = A^c(\alpha) [x^c-x^c_e(\alpha)] + B^c(\alpha) [y-r],\\[5pt]
u = C^c(\alpha) [x^c-x^c_e(\alpha)]+ D^c(\alpha) [y-r]+ u_e(\alpha), ~~~ \forall \alpha \in \Omega.
\end{array}
\end{equation}

Note that $r_e(\alpha)=y_e(\alpha)$, as a result $\delta y - \delta r=y-r$. The scheduling parameter $\alpha$ is treated as a parameter throughout the design process, and then it becomes a time-varying input signal to the gain-scheduled controller implementation through the dependence $\alpha(t)=p(y(t))$. Thus the gain-scheduled controller becomes
\begin{equation}\label{eqn_gs14}
\begin{array}{l}
\dot{x}^c = A^c(p(y)) [x^c-x^c_e(p(y))] + B^c(p(y)) [y-r],\\[5pt]
u = C^c(p(y)) [x^c-x^c_e(p(y))]+ D^c(p(y)) [y-r]+ u_e(p(y)).
\end{array}
\end{equation}

Linearization of (\ref{eqn_gs14}) about an equilibrium specified by $\alpha$ gives
\begin{equation}\label{eqn_gs15}
\begin{array}{l}
 \delta \dot{x}^c = A^c(\alpha) \delta x^c + B^c(\alpha) [y-r] \displaystyle -[A^c(\alpha) \frac{\partial x^c_e (\alpha)}{\partial \alpha}] \times [\frac{\partial p }{\partial y} (y_e(\alpha)) (y-r)], \\[5pt]
\displaystyle \delta u = C^c(\alpha) \delta x^c + D^c(\alpha) [y-r] \displaystyle  +[\frac{\partial u_e (\alpha)}{\partial \alpha} - C^c(\alpha) \frac{\partial x^c_e (\alpha)}{\partial \alpha}] \displaystyle \times [\frac{\partial p }{\partial y} (y_e(\alpha)) (y-r)].
\end{array}
\end{equation}

Comparing this to (\ref{eqn_gs9}), we see there are additional terms, which we refer to them as hidden coupling terms \cite{research-rugh-2000}.

In order to get rid of these terms we have to design the controller such that
\begin{equation}\label{eqn_gs16}
\begin{array}{l}
\displaystyle A^c(\alpha) \frac{\partial x^c_e (\alpha)}{\partial \alpha} = 0, \\[5pt]
\displaystyle \frac{\partial u_e (\alpha)}{\partial \alpha} - C^c(\alpha) \frac{\partial x^c_e (\alpha)}{\partial \alpha} =0.
\end{array}
\end{equation}

It is not always easy to come up with solutions to satisfy this condition. In order to make the design process easier, we can augment integrators at the plant input, so there is no need for equilibrium control value. By augmenting integrators at the plant (\ref{eqn_gs1}) input we obtain
\begin{equation} \label{eqn_gs61}
\begin{array}{l}
  \left[
       \begin{array}{c}
           \dot{x}^p(t) \\
           \dot{u}(t)
       \end{array}
    \right] =
     \left(
       \begin{array}{c}
           f^p(x^p(t),u(t)) \\
           -\eta_c u(t)
       \end{array}
    \right) +
     \left[
       \begin{array}{c}
           0 \\
           \eta_c \times I
           \end{array}
    \right] v(t).
\end{array}
\end{equation}

The new controller has the general form
\begin{equation}\label{eqn_gs110}
\begin{array}{c}
\dot{x}^c(t)= f^c(x^c(t),y(t), r(t)),\\[5pt]
v(t)=g^c(x^c(t),y(t), r(t)),
\end{array}
\end{equation}
with the input and output signals corresponding to the nonlinear plant (\ref{eqn_gs1}).

Combining (\ref{eqn_gs61}) and (\ref{eqn_gs110})leads to
\begin{equation} \label{eqn_gs111}
\begin{array}{l}
  \underbrace{\left[
       \begin{array}{c}
           \dot{x}^p(t) \\
           \dot{u}(t) \\
           \dot{x}^c(t)
       \end{array}
    \right]}_{\dot{x}} =
     \underbrace{\left(
       \begin{array}{c}
           f^p(x^p(t),u(t)) \\
           -\eta_c u(t)  \\
           f^c(x^c(t),g^p(x^p(t),u(t)) ,r(t))
       \end{array}
    \right)}_{f(x,r)} +
     \underbrace{\left[
       \begin{array}{c}
           0 \\
           \eta_c \times I \\
           0
       \end{array}
    \right]}_{B} v(t). \\[5pt]
v(t)=\underbrace{g^c(x^c(t),g^p(x^p(t),u(t)), r(t))}_{g(x,r)},
\end{array}
\end{equation}
Then the closed loop nonlinear system is

\begin{equation}\label{eqn_gs112}
\begin{array}{l}
\dot{x}(t)= f(x(t),r(t))+B g(x(t),r(t )),\\[5pt]
~~~~~~ =F(x(t),r(t))
\end{array}
\end{equation}

The augmented linear family of systems for (\ref{eqn_gs61}) becomes
\begin{equation}\label{eqn_gs62}
\begin{array}{l}
  \underbrace{\left[
       \begin{array}{c}
           \dot{x}^p(t) \\
           \dot{u}(t)
       \end{array}
    \right]}_{\dot{x}_{aug}} =
     \underbrace{\left[
       \begin{array}{cc}
           A^p(\alpha) &  B^p(\alpha)\\
              0      &  -\eta_c \times I
       \end{array}
    \right]}_{A_{aug}(\alpha)}
    \underbrace{\left[
       \begin{array}{c}
           \delta x^p \\
           \delta u
           \end{array}
    \right]}_{\delta x_{aug}}+
     \underbrace{\left[
       \begin{array}{c}
           0 \\
           \eta_c \times I
           \end{array}
    \right]}_{B_{aug}} v(t) ,\\[5pt]
\delta y = \underbrace{[C^p(\alpha), D^p(\alpha)]}_{C_{aug}(\alpha)}
   \underbrace{\left[
       \begin{array}{c}
           \delta x^p \\
           \delta u
           \end{array}
    \right]}_{\delta x_{aug}}.
\end{array}
\end{equation}
Now, the control realization for this system is
\begin{equation}\label{eqn_gs63}
\begin{array}{l}
\dot{x}^c = A^c(\alpha) x^c + B^c(\alpha) [y-r],\\[5pt]
v = C^c(\alpha) x^c+ D^c(\alpha) [y-r] ~~~ \forall \alpha \in \Omega.
\end{array}
\end{equation}
One of the options for control design is to set controller matrices as follows
\begin{equation}\label{eqn_gs17}
\begin{array}{ll}
A^c(\alpha)= A^c= -\epsilon_c I, & ~B^c(\alpha)=B^c=I, \\[5pt]
C^c(\alpha) = -K_i(\alpha), & ~D^c(\alpha)=-K_p(\alpha).
\end{array}
\end{equation}
which is a kind of PI control, where $K_i(\alpha)$ is the integral gain matrix, and $K_p(\alpha)$ is the proportional gain matrix.

Hence the control for the augmented system has the final form
\begin{equation} \label{eqn_gs64}
\begin{array}{l}
  \left[
       \begin{array}{c}
           \dot{x}^c \\
            v
       \end{array}
    \right] =
     \left[
       \begin{array}{ccc}
           -\epsilon_c I & I &  -I\\
           K_i(\alpha) & K_p(\alpha) & -K_p(\alpha)
       \end{array}
    \right] ~
     \left[
       \begin{array}{c}
            x^c \\
            y \\
            r
       \end{array}
    \right], ~~~ \forall \alpha \in \Omega.
\end{array}
\end{equation}

With these choices for control matrices, the control input is
\begin{equation}\label{eqn_gs19}
\begin{array}{l}
x^c =  \displaystyle \int \! \left(-\epsilon_c x^c +(y-r) \right) \, \mathrm{d}\tau,  \\[5pt]
v = - K_i(\alpha)  x^c - K_p(\alpha) (y-r), ~~~ \forall \alpha \in \Omega.
\end{array}
\end{equation}

\begin{figure}[!ht]
\centering
\includegraphics[width=0.5\textwidth]{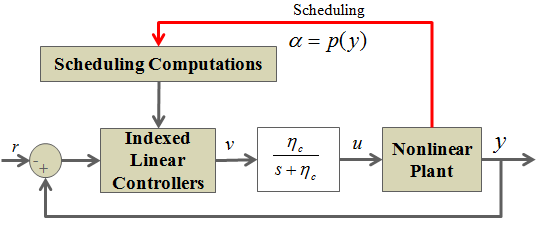}
\caption{Gain scheduling controller diagram}
\label{fig:GainSched_Cont_Struc}
\end{figure}

Figure~\ref{fig:GainSched_Cont_Struc}, shows schematically how the gain scheduling controller works.

The linearized closed loop system (\ref{eqn_gs62}) with controller (\ref{eqn_gs63}) becomes
%\end{multicols}
%\vspace{-21mm}
%\section*{{}}
%\begin{widetext}
\begin{equation} \label{eqn_gs21}
\begin{array}{l}
  \underbrace{\left[
       \begin{array}{c}
           \delta \dot{x}^p \\
           \delta \dot{u} \\
            \dot{x}^c
       \end{array}
    \right]}_{\dot{x}} =
     \underbrace{\left[
       \begin{array}{ccc}
           A^p(\alpha) &~~~ B^p(\alpha) &~~~ 0 \\
           \eta_c D^c(\alpha) C^p(\alpha) &~~~ -\eta_c I+ D^c(\alpha) D^p(\alpha) & ~~~ \eta_c C^c(\alpha) \\
           B^c(\alpha) C^p(\alpha)      &~~~ B^c(\alpha)D^p(\alpha)  &~~~ A^c(\alpha)
       \end{array}
    \right]}_{A_{cl}(\alpha)}
     \underbrace{\left[
       \begin{array}{c}
           \delta x^p \\
           \delta u \\
            x^c
       \end{array}
    \right]}_{x} +
    \underbrace{\left[
       \begin{array}{c}
            0 \\
           -\eta_c D^c(\alpha) \\
           - B^c(\alpha)
       \end{array}
    \right]}_{B_{cl}(\alpha)} \delta r, ~~~~~ \forall \alpha \in \Omega.
\end{array}
\end{equation}
%\end{widetext}

%\begin{multicols}{2}

For the case where we have simplified output $\delta y = \delta x^p$, (i.e. $C^p(\alpha)=I, D^p(\alpha)=0$) the linearized closed loop system (\ref{eqn_gs62}) with controller (\ref{eqn_gs64}) becomes
\begin{equation} \label{eqn_gs22}
\begin{array}{l}
  \underbrace{\left[
       \begin{array}{c}
           \delta \dot{x}^p \\
           \delta \dot{u} \\
            \dot{x}^c
       \end{array}
    \right]}_{\dot{x}} =
     \underbrace{\left[
       \begin{array}{ccc}
           A^p(\alpha) &~~~ B^p(\alpha) &~~~ 0 \\
           -\eta_c K_p(\alpha) &~~~ -\eta_c I & ~~~ -\eta_c K_i(\alpha) \\
           I      &~~~ 0  &~~~ -\epsilon_c I
       \end{array}
    \right]}_{A_{cl}(\alpha)}
     \underbrace{\left[
       \begin{array}{c}
           \delta x^p \\
           \delta u \\
            x^c
       \end{array}
    \right]}_{x} +\underbrace{\left[
       \begin{array}{c}
            0 \\
           \eta_c K_p(\alpha) \\
           - I
       \end{array}
    \right]}_{B_{cl}(\alpha)} \delta r, ~~~ \forall \alpha \in \Omega.
\end{array}
\end{equation}

%%%%%%%%%%%%%%%%%%%%%%%%%%%%%%%%%%%%%%%%%%%%%%
\subsection{Stability Analysis}

In closed loop system (\ref{eqn_gs22}), let $\delta r =0$, and consider the unforced linear time varying system
\begin{equation}\label{eqn_gs50}
\begin{array}{l}
\dot{x}=A_{cl}(\alpha) x, ~~ x(0)= x_0, ~~ \forall \alpha \in \Omega, \\[5pt]
\delta y= \delta x^p.
\end{array}
\end{equation}

\newtheorem{ass}{Assumption}
\begin{ass} \label{ass1}
The matrix $A_{cl}$ is bounded and Lipschitz continuous as follows
\begin{equation}\label{eqn_gs51}
\begin{array}{l}
||A_{cl}(t)|| \leq k_A, ~~ \forall t>0, \\[5pt]
||A_{cl}(t)-A_{cl}(\tau)|| \leq L_A ||t- \tau||, ~~ \forall t , \tau >0,
\end{array}
\end{equation}
\end{ass}

\begin{ass} \label{ass2}
The constant eigenvalues of matrix $A_{cl}(y)$ are uniformly bounded away from the closed complex right-half plane for all constant $y$.
\end{ass}

\newtheorem{thm}{Theorem}
\begin{thm} \label{thm1}
Consider system (\ref{eqn_gs50}), under assumptions \ref{ass1} and \ref{ass2}, then there exists constants $m$, $\lambda$, and $\epsilon > 0$ such that if
\begin{equation}\label{eqn_gs52}
||\dot{y}(t)|| \leq \epsilon_y, ~~ \forall t \in [0,T],
\end{equation}
then
\begin{equation}\label{eqn_gs53}
 ||x(t)|| \leq me^{-\lambda t} ||x_0||,  ~~\forall t \in [0,T].
\end{equation}
\end{thm}

To analyse the stability of the nonlinear closed-loop system, we use a technique known as "global linearization" developed in \cite{lmi-boyd-1994}.

\begin{thm}\label{thm2}
Consider nonlinear system (\ref{eqn_gs112}), and assume there are a family of equilibrium points $(x_e,r_e)$ such that $F(x_e,r_e)=0$. Then $A^{nl}_{cl} = \frac{\partial F}{\partial x} \in S, ~\forall x$, where $S$ is a polytope, and it is described by a list of its vertices, i.e. in the form
\begin{equation}\label{eqn_gs113}
 S:= \textbf{Co}\{ A^{nl}_{cl_1}, ..., A^{nl}_{cl_L} \},
\end{equation}
where $A^{nl}_{cl_i}$s are obtained by linearizing nonlinear system (\ref{eqn_gs112}) near equilibrium points (steady state condition), and also non-equilibrium points (transient condition). Now, assume their exist a common symmetric positive definite matrix $P=P^\mathsf{T} > 0$ such that: \\
\begin{equation}\label{eqn_gs54}
 P A^{nl}_{cl_i}+A^{nl \mathsf{T}}_{cl_i} P < 0, ~~~ \forall i \in \{ 1,2,..., L \}.
\end{equation}
then system (\ref{eqn_gs112}) is absolutely stable. Since by design $A_{cl}(\alpha) \in S, ~\forall \alpha$, then system (\ref{eqn_gs50}) is also stable.
\end{thm}

%\begin{proof}
%A Lyapunov candidate function chosen as
%\begin{equation}\label{eqn_gs154}
% V= x^\mathsf{T} P  x
%\end{equation}
%where its time-derivative is given by
%\begin{equation}\label{eqn_gs155}
%\begin{array}{l}
%\dot{V}= \dot{x}^\mathsf{T} P x + x^\mathsf{T} P \dot{x} \\[5pt]
%~~~ = x^\mathsf{T} \left( A_{cl}^\mathsf{T}(\alpha) P + P A_{cl}(\alpha) \right) x,
%\end{array}
%\end{equation}
%using equation (\ref{eqn_gs55})
%\begin{equation}\label{eqn_gs156}
%\dot{V} < \|x\|^2,
%\end{equation}
%So the system (\ref{eqn_gs50}) is asymptotically stable.
%\hfill\(\Box\)
%\end{proof}

%\newtheorem{rmk}{Remark}
%\begin{rmk} \label{rmk2}
%Since in inequality (\ref{eqn_gs54}), the $A_{cl_i}$'s we are using contains a wide range of linearizations of the nonlinear system (\ref{eqn_gs1}), including equilibrium and non-equilibrium linearizations, we can claim that the closed loop system (\ref{eqn_gs61}) and (\ref{eqn_gs19}) is also stable.
%\end{rmk}

If there are no hidden coupling terms involving $\delta y$, then the design of a stabilizing linear controller family can be assumed to guarantee stability of the linearized closed-loop system in a neighborhood of every $\alpha \in \Omega$. The closed-loop system is not restricted to remain in a neighborhood of any single equilibrium, but is assumed to be \emph{slowly-varying} and to have initial state sufficiently close to some equilibrium in S. Then the conclusion is that the closed-loop system remains in a neighborhood of the equilibrium manifold \cite{research-rugh-2000}. Using results developed in \cite{gainsched-shamma-1988}, we can figure out if a system is slowly-varying or not. Here we rewrite theorem Theorem 12 from \cite{research-rugh-2000}:

\begin{thm} \label{thm3}
For plant (\ref{eqn_gs1}), suppose the gain-scheduled controller (\ref{eqn_gs14}) is such that there are no hidden coupling terms and the eigenvalues of the linearized closed-loop system satisfy $Re[\lambda] \leq - \epsilon < 0$ for every $\alpha \in \Omega$. Then given $\rho > 0$ there exist positive constants $\mu$ and $\gamma$ such that the response of the nonlinear closed-loop system satisfies the following property. If the exogenous signal $|| \dot{r}(t)|| < \mu$, for $t \geq 0$, and if for some $\alpha \in \Omega$,
\begin{equation} \label{eqn_gs23}
\begin{array}{l}
 \left|\left|
 \left[
       \begin{array}{c}
           x^p(0) \\
           u(0) \\
           x^c(0)
       \end{array}
    \right] -
     \left[
       \begin{array}{c}
           x^p_e(\alpha) \\
           u_e(\alpha) \\
           0
       \end{array}
    \right]
 \right|\right| < \gamma,
\end{array}
\end{equation}
 then
 \begin{equation} \label{eqn_gs24}
\begin{array}{l}
  \left|\left|
   \left[
       \begin{array}{c}
           x^p(t) \\
           u(t) \\
           x^c(t)
       \end{array}
    \right] -
     \left[
       \begin{array}{c}
           x^p_e(p(y(t))) \\
           u_e(p(y(t))) \\
           0
       \end{array}
    \right]
    \right|\right| < \rho,   ~~\forall t \geq 0.
\end{array}
\end{equation}
\end{thm}
%%%%%%%%%%%%%%%%%%%%%%%%%%%%%%%%%%%%%%%%%%%%%%%%%%
\subsection{Integration Anti-Windup}

It is desirable to have integral action in the controller since the presence of an integral term eliminates steady state error in the controlled variable. Since the allowable values for control inputs are limited, if any controller reaches its limit, and error is produced form the difference of the control signal and the actual limited signal applied to the plant. This phenomenon is known as integral wind-up. Because of this, Integral Wind-Up Protection (IWUP) is used to reduce the effect of the integral term of the controller.

An approach to IWUP from \cite{development-martin-2008, control-csnak-2010} was adopted for our control problem. The main idea with this approach is to decrease the error seen by the integrators. This allows the integrator to increase to an appropriate value and decrease the size of the instantaneous change in magnitude when the controller becomes saturated. First, the generated control signal to the fuel-metering valve is subtracted from the saturation value. The resulting difference is then amplified by an integral feedback gain (IFB) and subtracted from the input to the integrator. The IFB is empirically tuned to provide adequate performance. The IFB is not gain scheduled, a constant value is sufficient for good performance.

%%%%%%%%%%%%%%%%%%%%%%%%%%%%%%%%%%%%%%%%%%%%%%%%%%%%%%%%%%%%%%%%%%%%%%
\section{Turboshaft Engine Example}

We apply the developed gain-scheduling controller to a physics-based model of a turboshaft engine driving a variable pitch propeller developed in \cite{fitzgerald-model-2012, pakmehr-decentmodel-2011}. For a standard day at sea level condition we found five equilibrium points for linearizing the dynamics near them. The linearization matrices for these five equilibrium points and steady state values of the engine variables and control parameters are:

\begin{itemize}

	\item Equilibrium Point 1 (Full Thrust):\\
	$u_1^*=1.0,~ u_2^*=16 ~(\text{deg}),~ x_1^*=1.0,~ x_2^*=0.9524,~ T^*=255.8685 ~(N),~ \alpha^*=1.3810,$ and the matrices are
	\begin{eqnarray} \label{eqn_gs70}
%\begin{split}
\begin{array}{c}
  A_1=
   \left[
       \begin{array}{cc}
           -5 & 0 \\
            3.5 & -2.3
       \end{array}
    \right],~ B_1=
    \left[
       \begin{array}{cc}
           1.4 & 0 \\
           0.63  & -0.085
       \end{array}
       \right], ~ C_1=I, \\[10pt]
   Ki_1=
   \left[
       \begin{array}{cc}
           0.7 & 0.7 \\
           0.7 & 0.6
       \end{array}
    \right],~ Kp_1=
    \left[
       \begin{array}{cc}
           1.2 & 1.2 \\
           1.2 & 1.2
       \end{array}
       \right].
\end{array}
%\end{split}
\end{eqnarray}

\item Equilibrium Point 2:\\
	$u_1^*=0.7,~ u_2^*=16 ~(\text{deg}),~ x_1^*=0.8826,~ x_2^*=0.6263,~ T^*=181.9711 ~(N),~ \alpha^*=1.0822,$ and the matrices are
	\begin{eqnarray} \label{eqn_gs71}
%\begin{split}
\begin{array}{c}
  A_2=
   \left[
       \begin{array}{cc}
           -2.83 & -0.0008 \\
            1.20 & -2.10
       \end{array}
    \right],~ B_2=
    \left[
       \begin{array}{cc}
           1.14  &      0 \\
           0.78  & -0.054
       \end{array}
       \right], ~ C_2=I, \\[10pt]
   Ki_2=
   \left[
       \begin{array}{cc}
           0.6 & 0.6 \\
           0.6 & 0.5
       \end{array}
    \right],~ Kp_2=
    \left[
       \begin{array}{cc}
           1.1 & 1.1 \\
           1.1 & 1.1
       \end{array}
       \right].
\end{array}
%\end{split}
\end{eqnarray}

\item Equilibrium Point 3 (Cruise):\\
	$u_1^*=0.4685,~ u_2^*=16 ~(\text{deg}),~ x_1^*=0.7264,~ x_2^*=0.5,~ T^*=70.5125 ~(N),~ \alpha^*=0.8818,$ and the matrices are
	\begin{eqnarray} \label{eqn_gs72}
%\begin{split}
\begin{array}{c}
  A_3=
   \left[
       \begin{array}{cc}
           -1.9 & 0.061 \\
            0.45 & -1.1
       \end{array}
    \right],~ B_3=
    \left[
       \begin{array}{cc}
           1.57 & 0 \\
           0.3 & -0.023
       \end{array}
       \right], ~ C_3=I, \\[10pt]
   Ki_3=
   \left[
       \begin{array}{cc}
           0.5 & 0.5 \\
           0.5 & 0.4
       \end{array}
    \right],~ Kp_3=
    \left[
       \begin{array}{cc}
           1 & 1 \\
           1 & 1
       \end{array}
       \right].
\end{array}
%\end{split}
\end{eqnarray}

\item Equilibrium Point 4:\\
	$u_1^*=0.3,~ u_2^*=16 ~(\text{deg}),~ x_1^*=0.5327,~ x_2^*=0.3678,~ T^*=38.155 ~(N),~ \alpha^*=0.6473,$ and the matrices are
	\begin{eqnarray} \label{eqn_gs73}
%\begin{split}
\begin{array}{c}
  A_4=
   \left[
       \begin{array}{cc}
           -0.85 & 0.032 \\
            0.32 & -0.64
       \end{array}
    \right],~ B_4=
    \left[
       \begin{array}{cc}
           1.1 & 0 \\
           0.17 & -0.011
       \end{array}
       \right], ~ C_4=I, \\[10pt]
   Ki_4=
   \left[
       \begin{array}{cc}
           0.4 & 0.4 \\
           0.4 & 0.3
       \end{array}
    \right],~ Kp_4=
    \left[
       \begin{array}{cc}
           0.8 & 0.8 \\
           0.8 & 0.8
       \end{array}
       \right].
\end{array}
%\end{split}
\end{eqnarray}

\item Equilibrium Point 5 (Idle):\\
	$u_1^*=0.145,~ u_2^*=16 ~(\text{deg}),~ x_1^*=0.295,~ x_2^*=0.161,~ T^*=7.317 ~(N),~ \alpha^*=0.3361,$ and the matrices are
	\begin{eqnarray} \label{eqn_gs74}
%\begin{split}
\begin{array}{c}
  A_5=
   \left[
       \begin{array}{cc}
           -0.38 & -0.0008 \\
            0.26  & -0.34
       \end{array}
    \right], ~ B_5=
    \left[
       \begin{array}{cc}
           0.7 & 0 \\
           0.1 & -0.0024
       \end{array}
       \right], ~C_5= I, \\[10pt]
   Ki_5=
   \left[
       \begin{array}{cc}
           0.3 & 0.3 \\
           0.3 & 0.2
       \end{array}
    \right],~ Kp_5=
    \left[
       \begin{array}{cc}
           0.6 & 0.6 \\
           0.6 & 0.6
       \end{array}
       \right].
\end{array}
%\end{split}
\end{eqnarray}
	
\end{itemize}

Other controller parameters are set to
\begin{equation}\label{eqn_gs75}
 \epsilon_c=1,  ~~ \eta_c=3, ~~ Q= 3 \times I.
\end{equation}

To show the stability of the closed loop system, 20 different (10 equilibrium, and 10 non-equilibrium) linearizations have been used, to solve inequality (\ref{eqn_gs54}), in Matlab with the aid of YALMIP \cite{YALMIP-lofberg-2004} and SeDuMi \cite{sedumi-Sturm-2001} packages. The numerical value for the common matrix $P$ is:

\begin{eqnarray} \label{eqn_gs76}
P = \left[
\begin{array}{cccccc}
    0.639 &  0.035 &  0.121 & -0.015 & -0.073 & -0.036  \\
    0.034 &  0.391 &  0.036 & -0.002 & -0.103 & -0.029  \\
    0.121 &  0.036 &  0.184 & -0.048 & -0.029 & -0.017  \\
   -0.015 & -0.002 & -0.048 &  0.130 &  0.028 &  0.022  \\
   -0.073 & -0.103 & -0.029 &  0.028 &  0.322 &  0.028  \\
   -0.036 & -0.029 & -0.017 &  0.022 &  0.028 &  0.298
\end{array}
\right]
\end{eqnarray}

\begin{figure}[!ht]
\centering
\begin{minipage}[l]{3.2in}
\centering
\resizebox{3.2in}{!}{\includegraphics{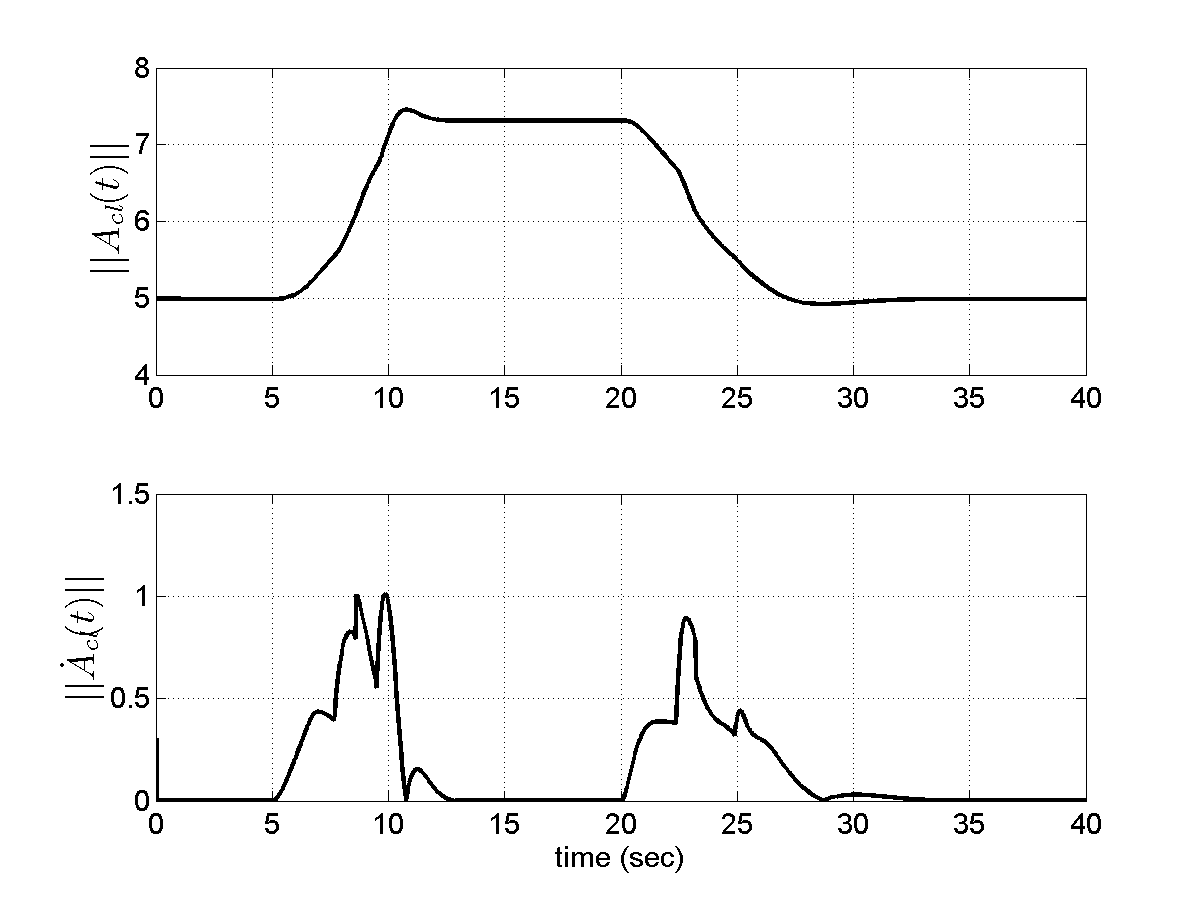}}
\caption{Norm of closed-loop system matrix ($||A_{cl}(t)||$), and its rate of change ($||\dot{A}_{cl}(t)||$) }\label{fig_cg_01}
\end{minipage}
\begin{minipage}[r]{3.2in}
\centering
\resizebox{3.2in}{!}{\includegraphics{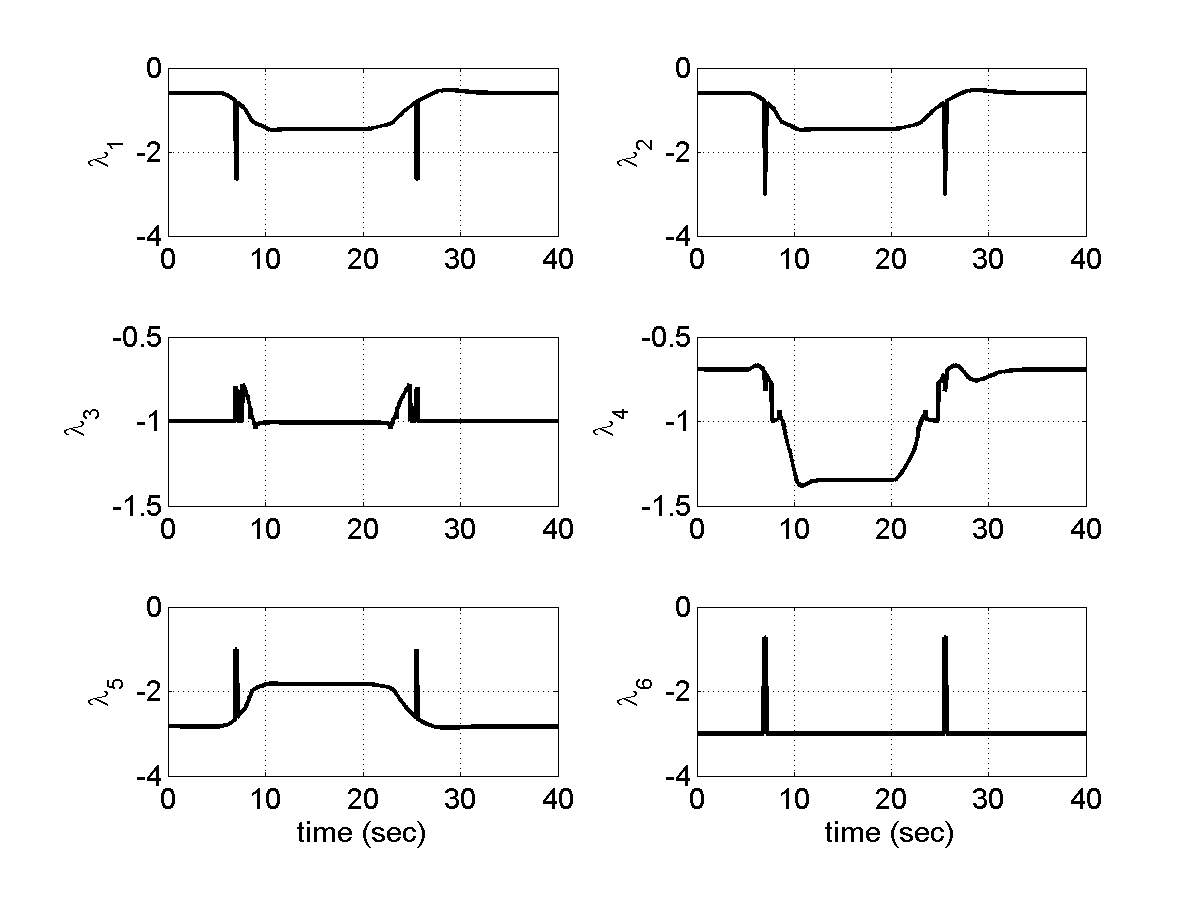}}
\caption{Closed-loop system eigenvalues ($\lambda[A_{cl}(\alpha)]$)}\label{fig_cg_02}
\end{minipage}
\end{figure}

\begin{figure}[!ht]
\centering
\begin{minipage}[l]{3.2in}
\centering
\resizebox{3.2in}{!}{\includegraphics{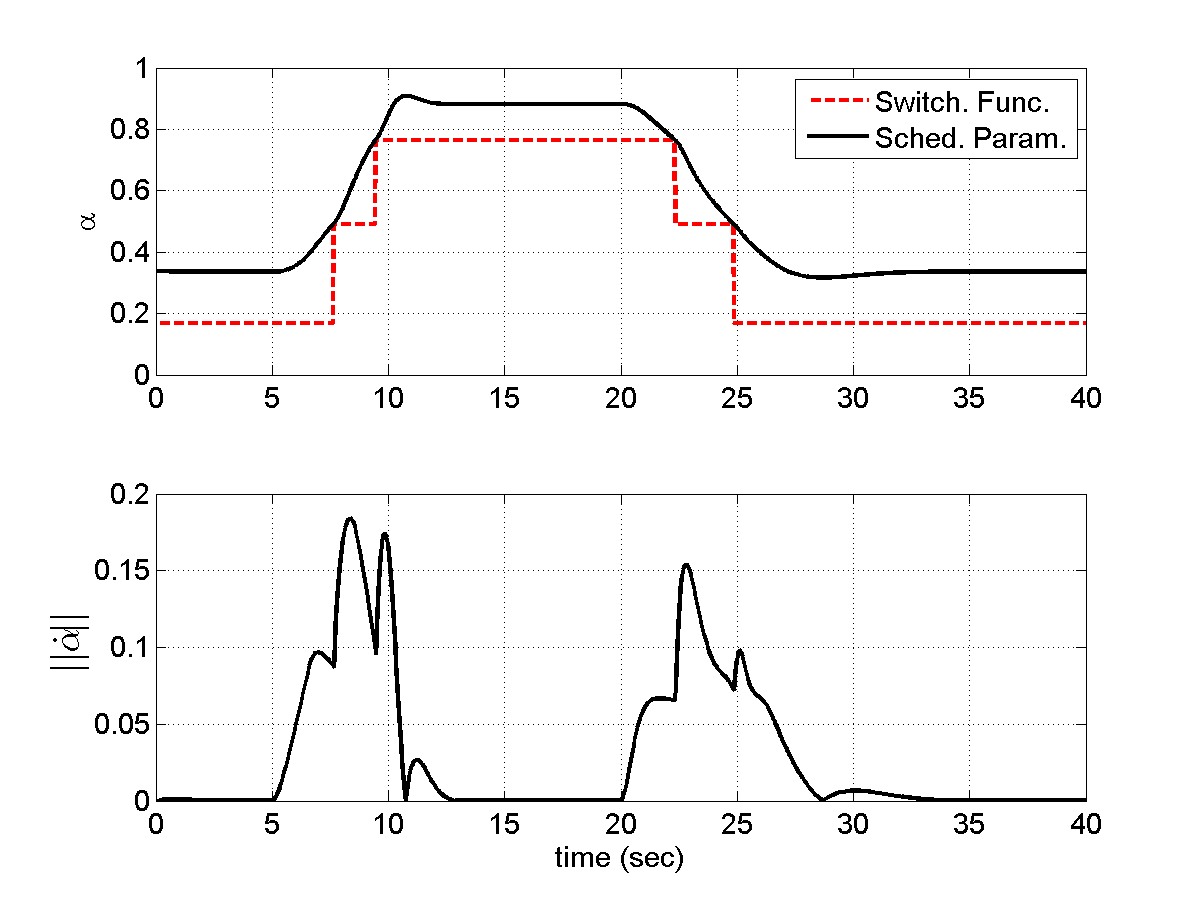}}
\caption{Scheduling Parameter ($\alpha(t)=||x(t)||$) and its rate of change ($\dot{\alpha}(t)=\frac{x^T \dot{x}}{||x(t)||}$)}\label{fig_cg_03}
\end{minipage}
\begin{minipage}[r]{3.2in}
\centering
\resizebox{3.2in}{!}{\includegraphics{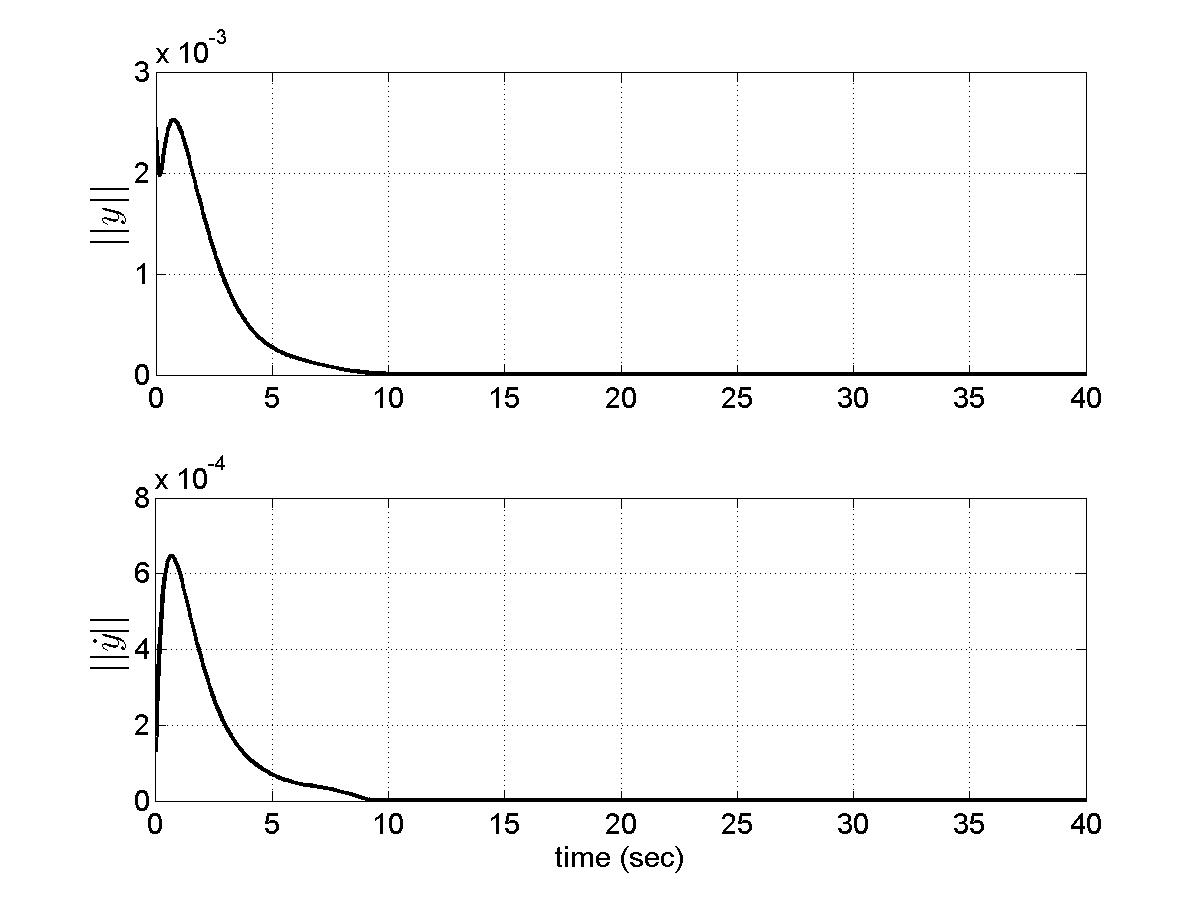}}
\caption{Norm of measured output of the system ($||y(t)||$), and its rate of change ($||\dot{y}(t)||$) }\label{fig_cg_04}
\end{minipage}
\end{figure}

\begin{figure}[!ht]
\centering
\begin{minipage}[l]{3.2in}
\centering
\resizebox{3.2in}{!}{\includegraphics{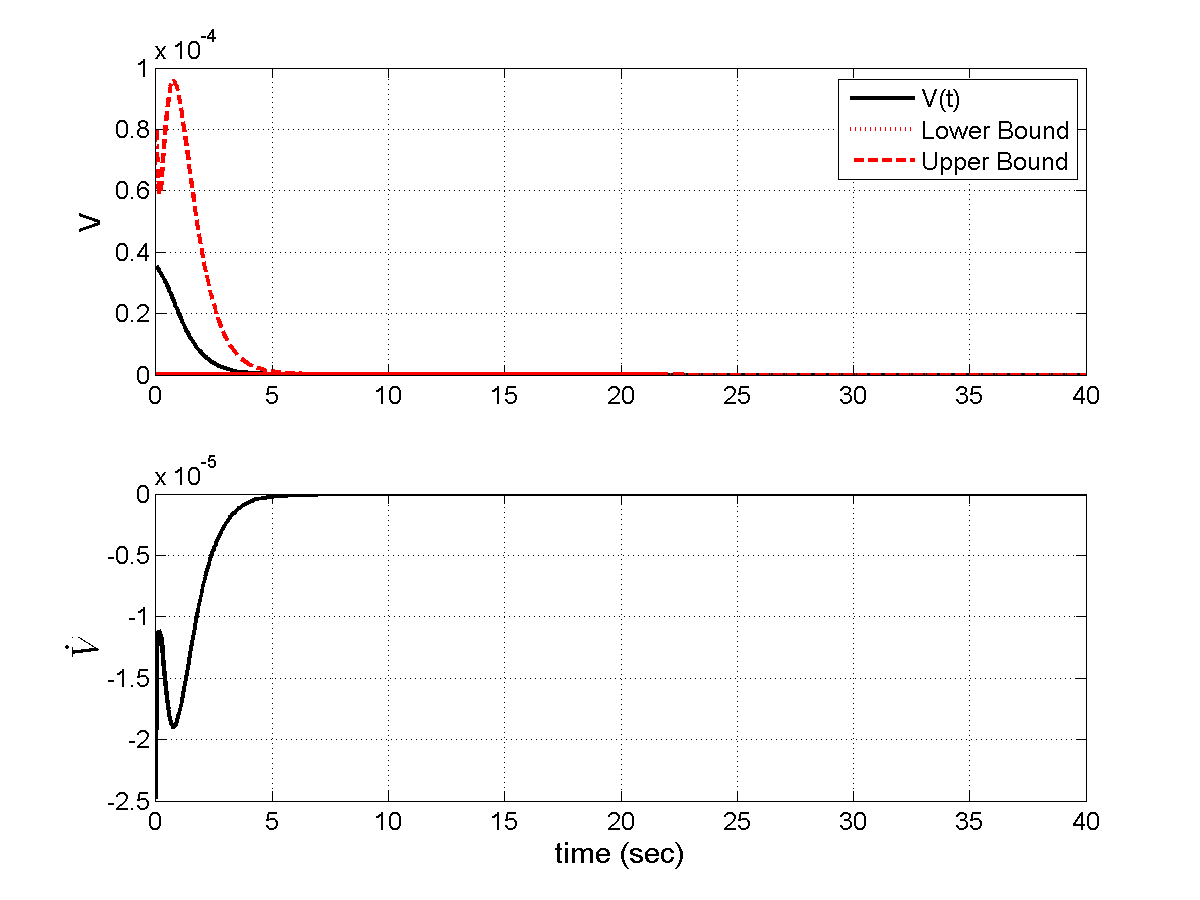}}
\caption{History of the unforced closed-loop system Lyapunov function $V(t)$, and its rate of change $\dot{V}(t)$ }\label{fig_cg_08}
\end{minipage}
\begin{minipage}[r]{3.2in}
\centering
\resizebox{3.2in}{!}{\includegraphics{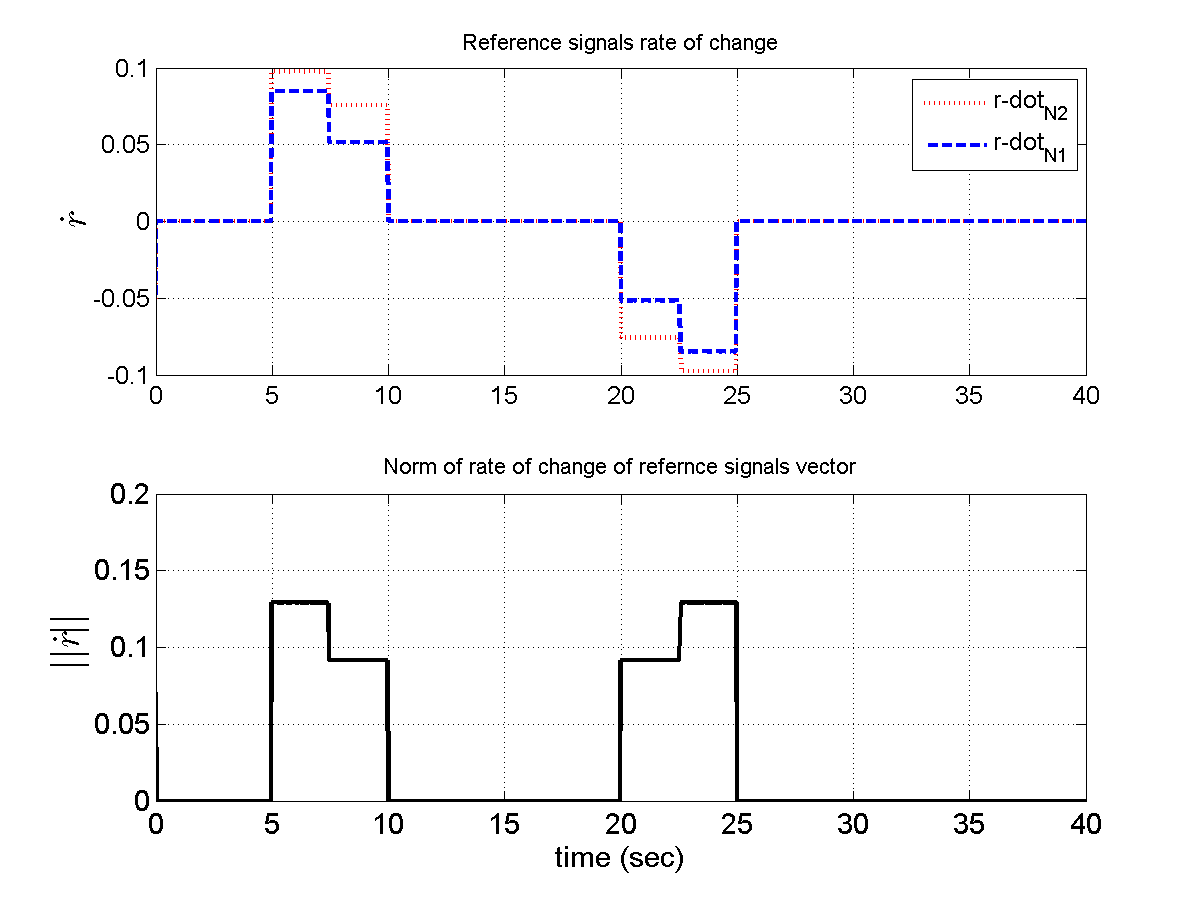}}
\caption{Rate of change of reference signals ($\dot{r}$)}\label{fig_cg_09}
\end{minipage}
\end{figure}

\begin{figure}[!ht]
\centering
\begin{minipage}[l]{3.2in}
\centering
\resizebox{3.2in}{!}{\includegraphics{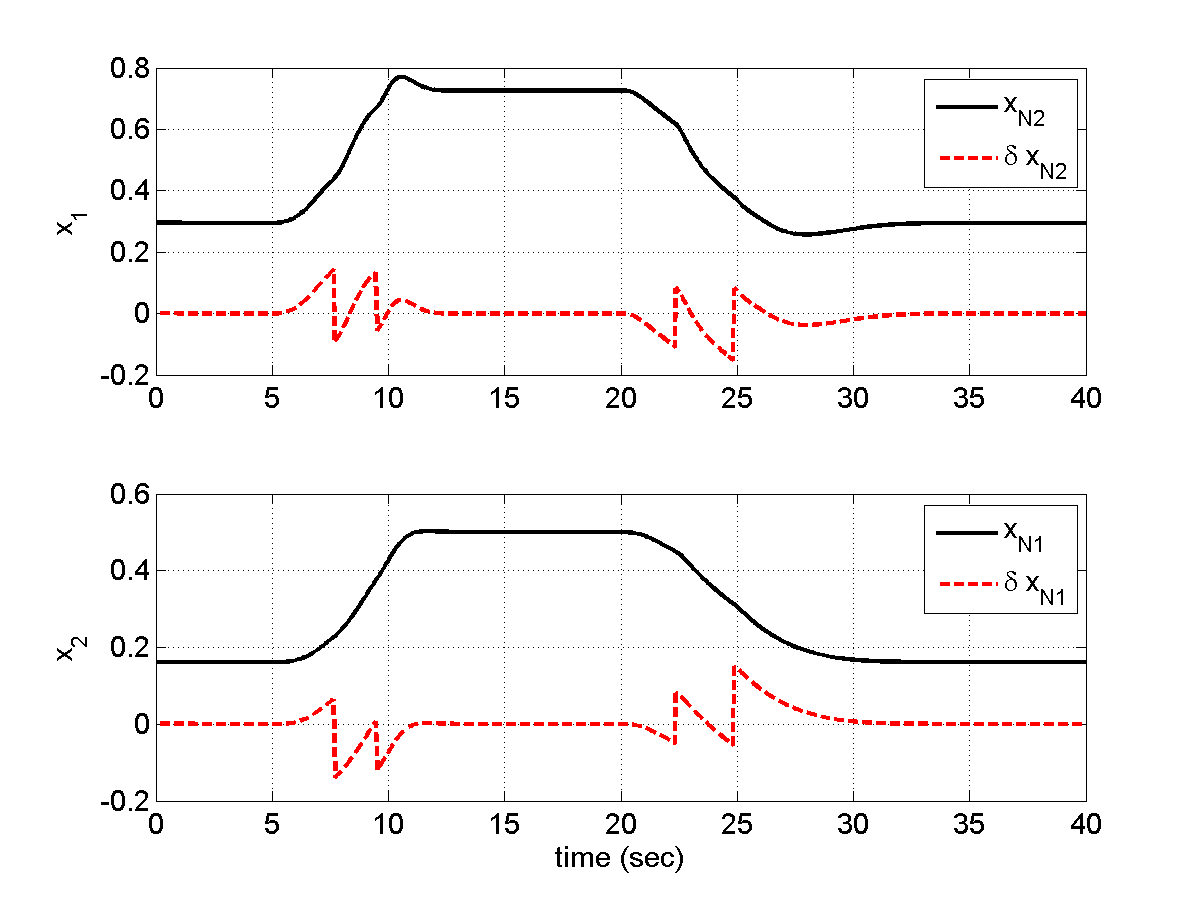}}
\caption{Plant states: core and fan spool speeds} \label{fig_cg_10}
\end{minipage}
\begin{minipage}[r]{3.2in}
\centering
\resizebox{3.2in}{!}{\includegraphics{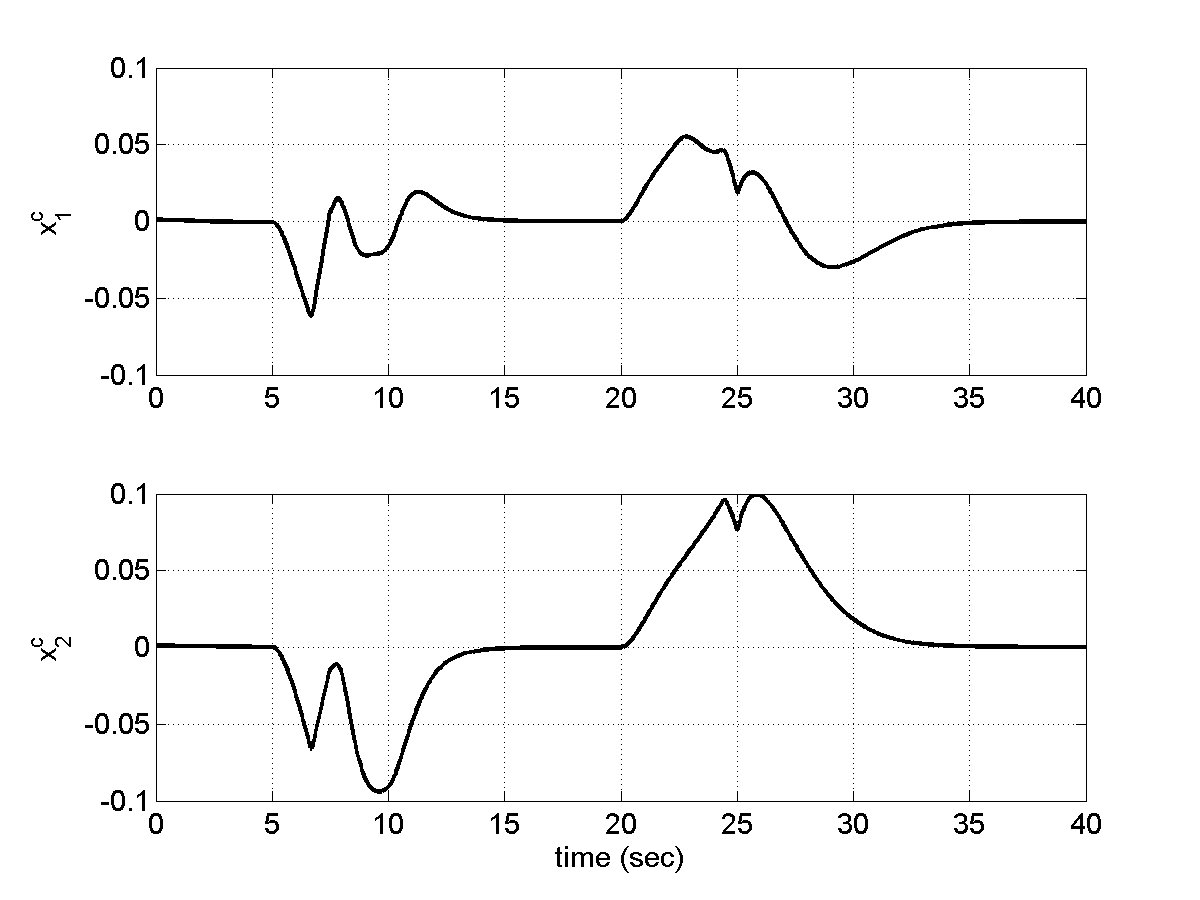}}
\caption{Controller states}\label{fig_cg_11}
\end{minipage}
\end{figure}

\begin{figure}[!ht]
\centering
\begin{minipage}[l]{3.2in}
\centering
\resizebox{3.2in}{!}{\includegraphics{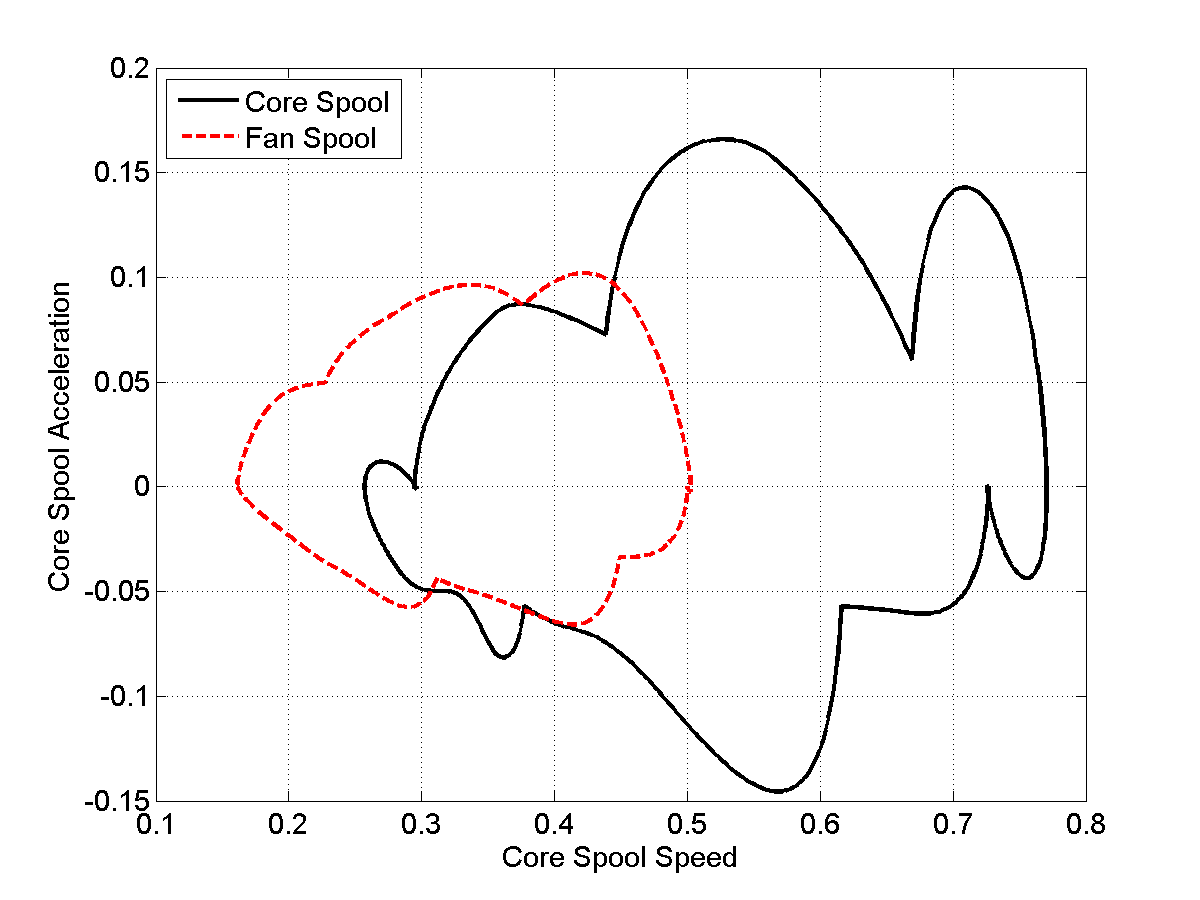}}
\caption{Core and fan spool speeds vs. core and fan spool accelerations}\label{fig_cg_12}
\end{minipage}
\begin{minipage}[r]{3.2in}
\centering
\resizebox{3.2in}{!}{\includegraphics{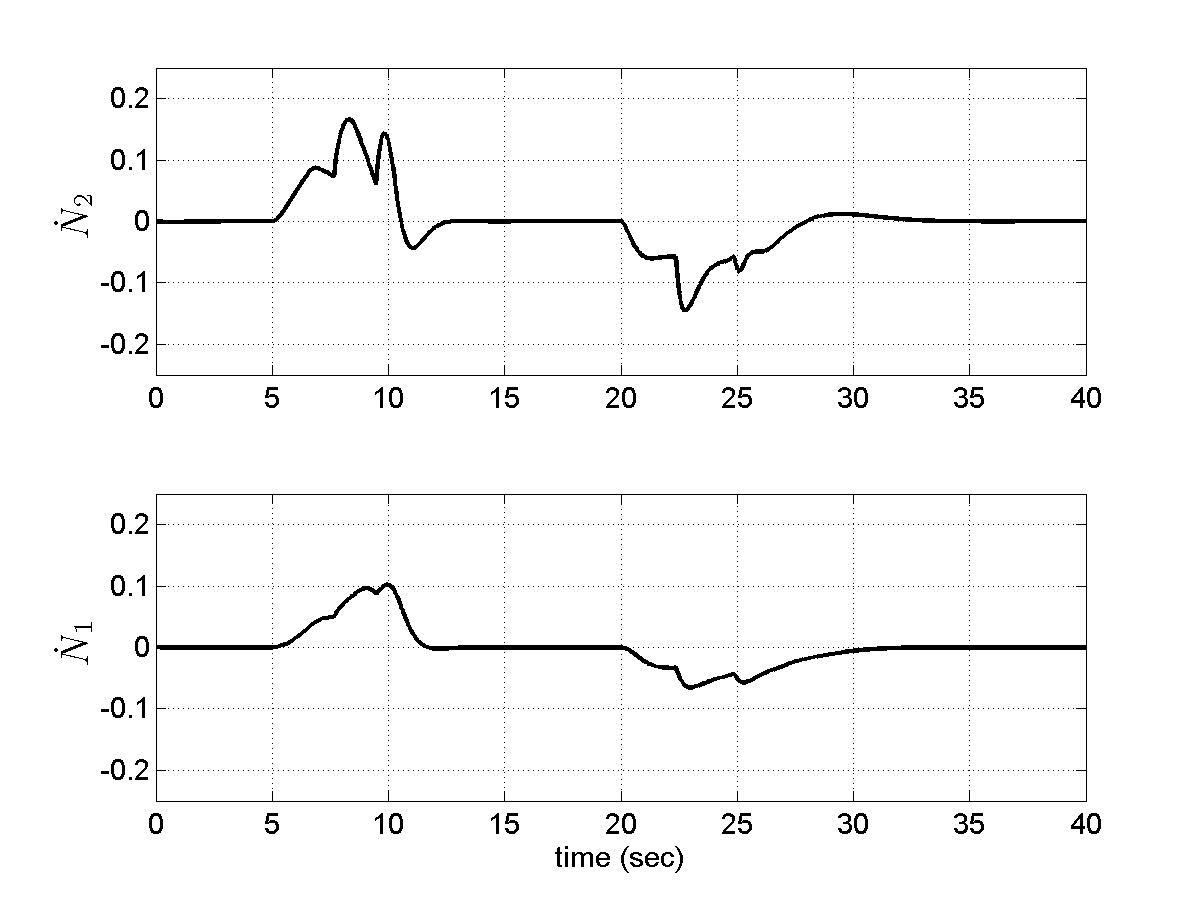}}
\caption{Core and fan spools accelerations}\label{fig_cg_13}
\end{minipage}
\end{figure}

\begin{figure}[!ht]
\centering
\begin{minipage}[l]{3.2in}
\centering
\resizebox{3.2in}{!}{\includegraphics{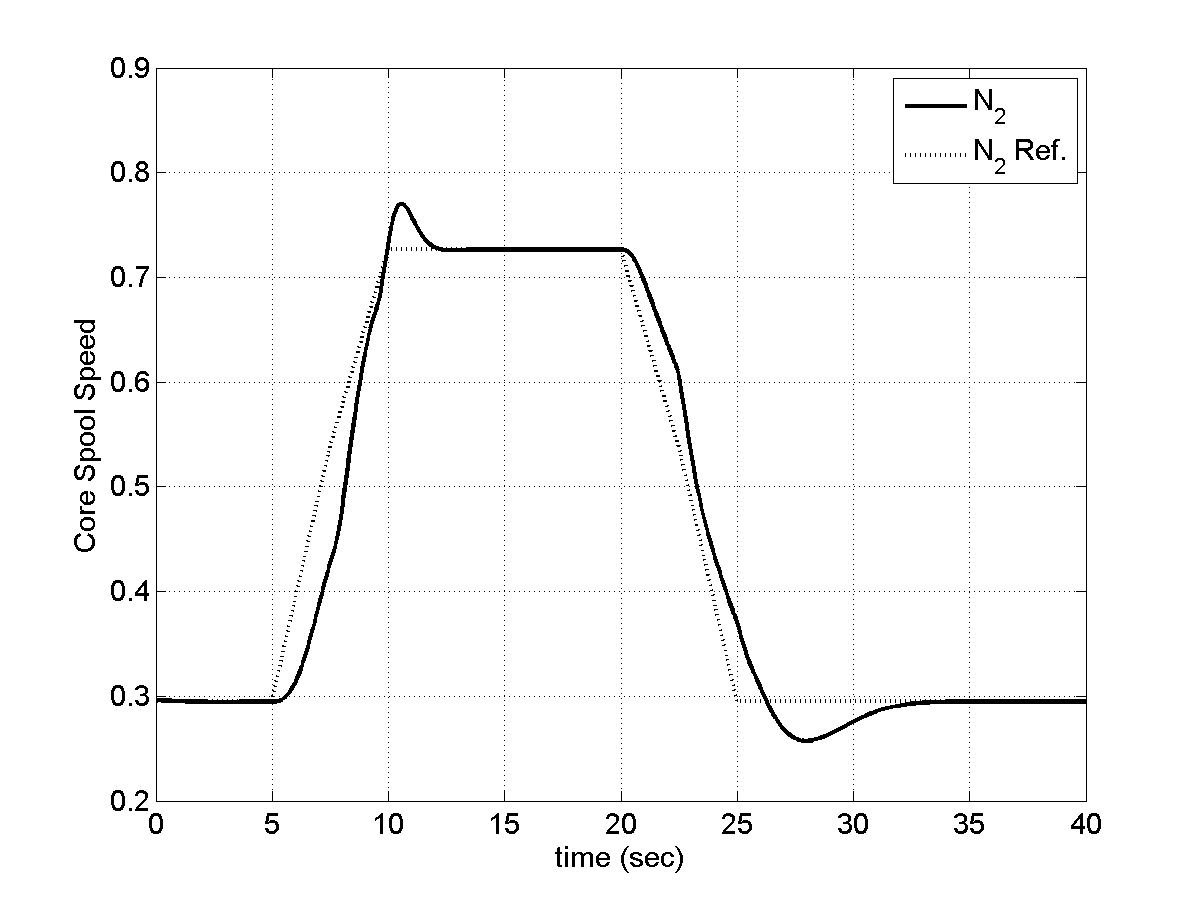}}
\caption{Output: core spool speed and its reference signal}\label{fig_cg_14}
\end{minipage}
\begin{minipage}[r]{3.2in}
\centering
\resizebox{3.2in}{!}{\includegraphics{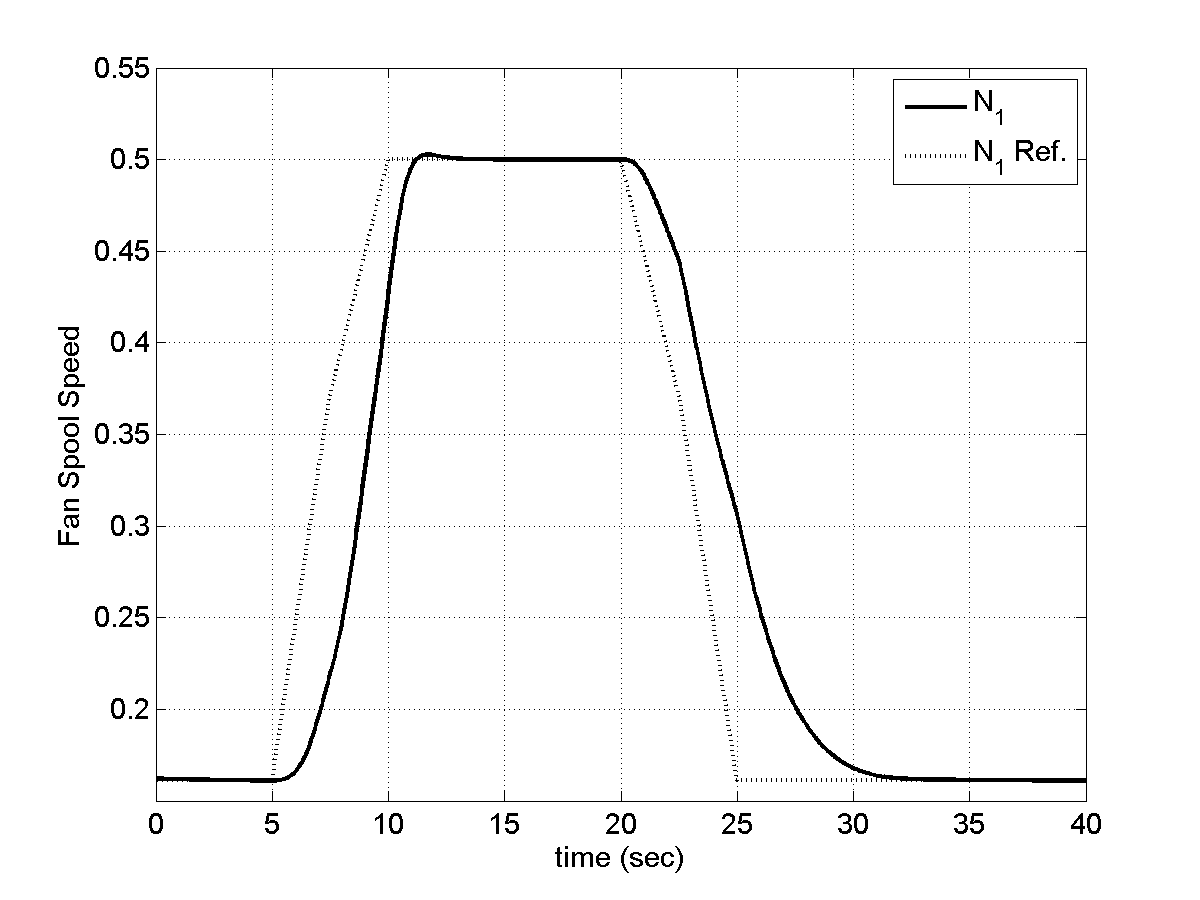}}
\caption{Output: fan spool speed and its reference signal}\label{fig_cg_15}
\end{minipage}
\end{figure}

\begin{figure}[!ht]
\centering
\begin{minipage}[l]{3.2in}
\centering
\resizebox{3.2in}{!}{\includegraphics{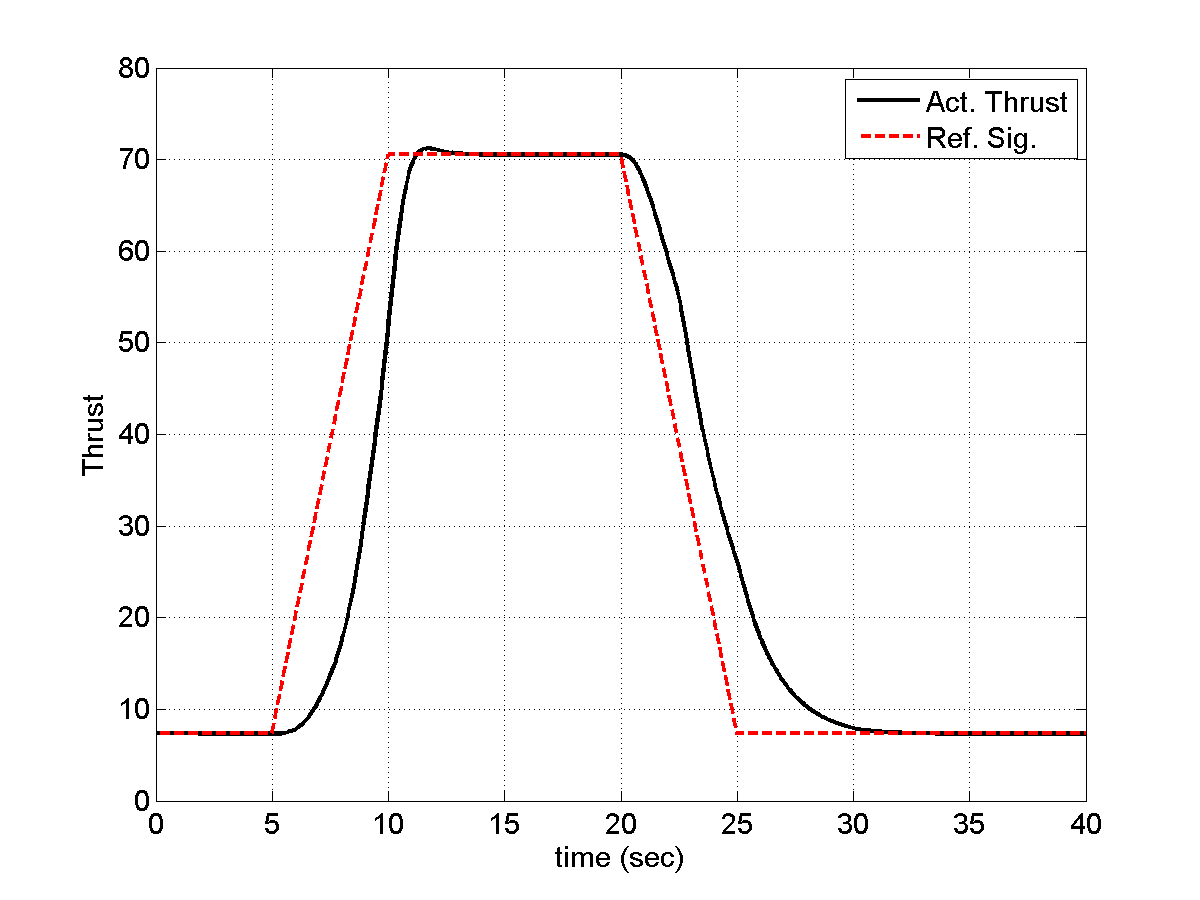}}
\caption{Thrust and its reference signal}\label{fig_cg_16}
\end{minipage}
\begin{minipage}[r]{3.2in}
\centering
\resizebox{3.2in}{!}{\includegraphics{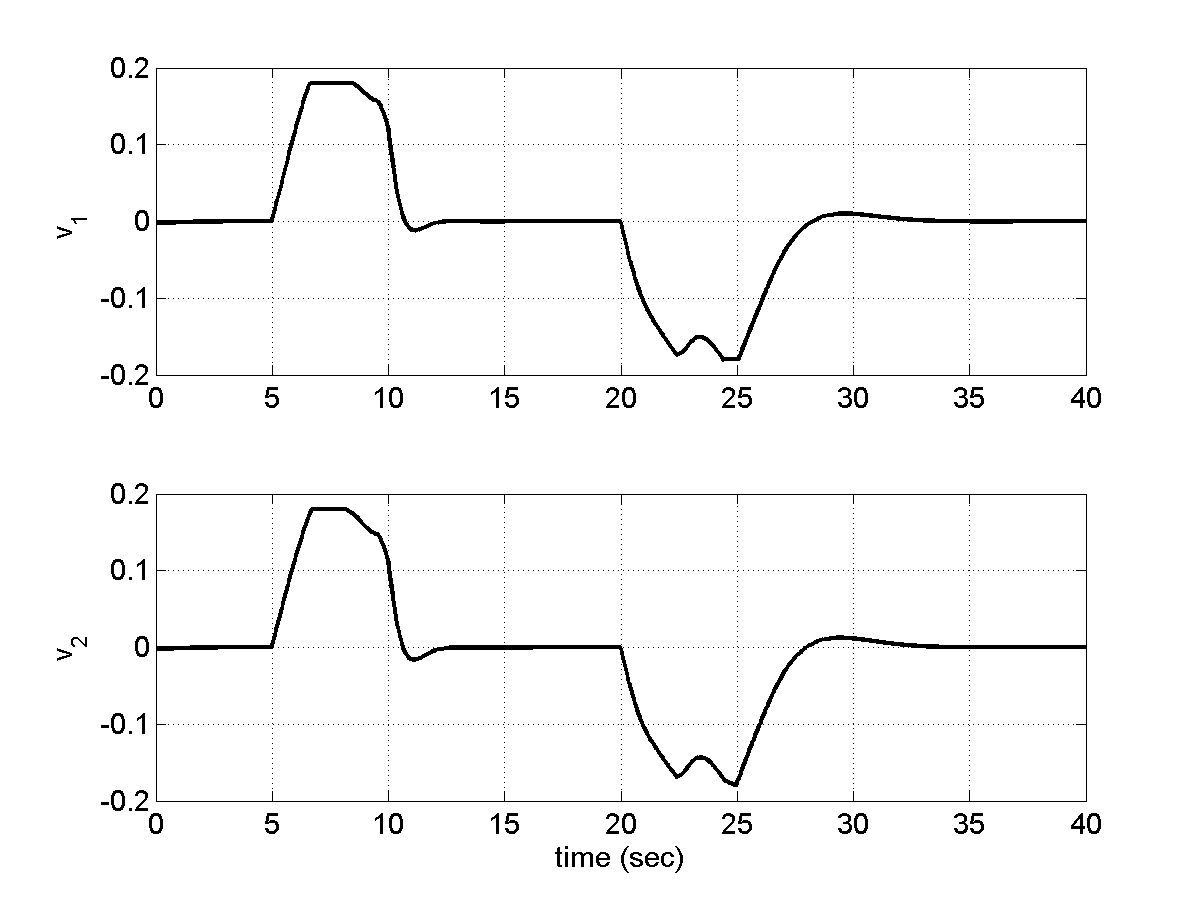}}
\caption{Control inputs to the augmented system ($v(t)$)}\label{fig_cg_17}
\end{minipage}
\end{figure}

\begin{figure}[!ht]
\centering
\begin{minipage}[l]{3.2in}
\centering
\resizebox{3.2in}{!}{\includegraphics{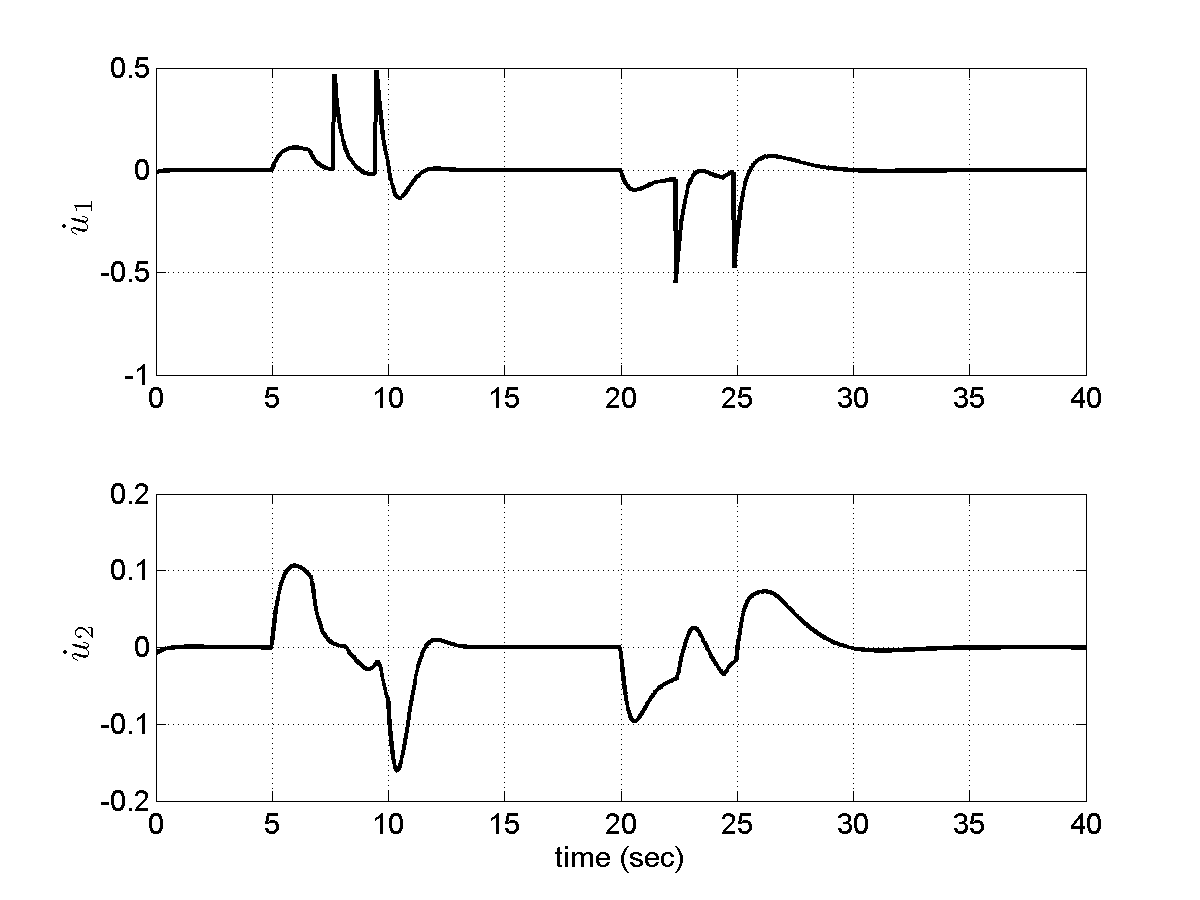}}
\caption{Rate of change for fuel and prop pitch angle control inputs ($\dot{u}(t)$)}\label{fig_cg_18}
\end{minipage}
\begin{minipage}[r]{3.2in}
\centering
\resizebox{3.2in}{!}{\includegraphics{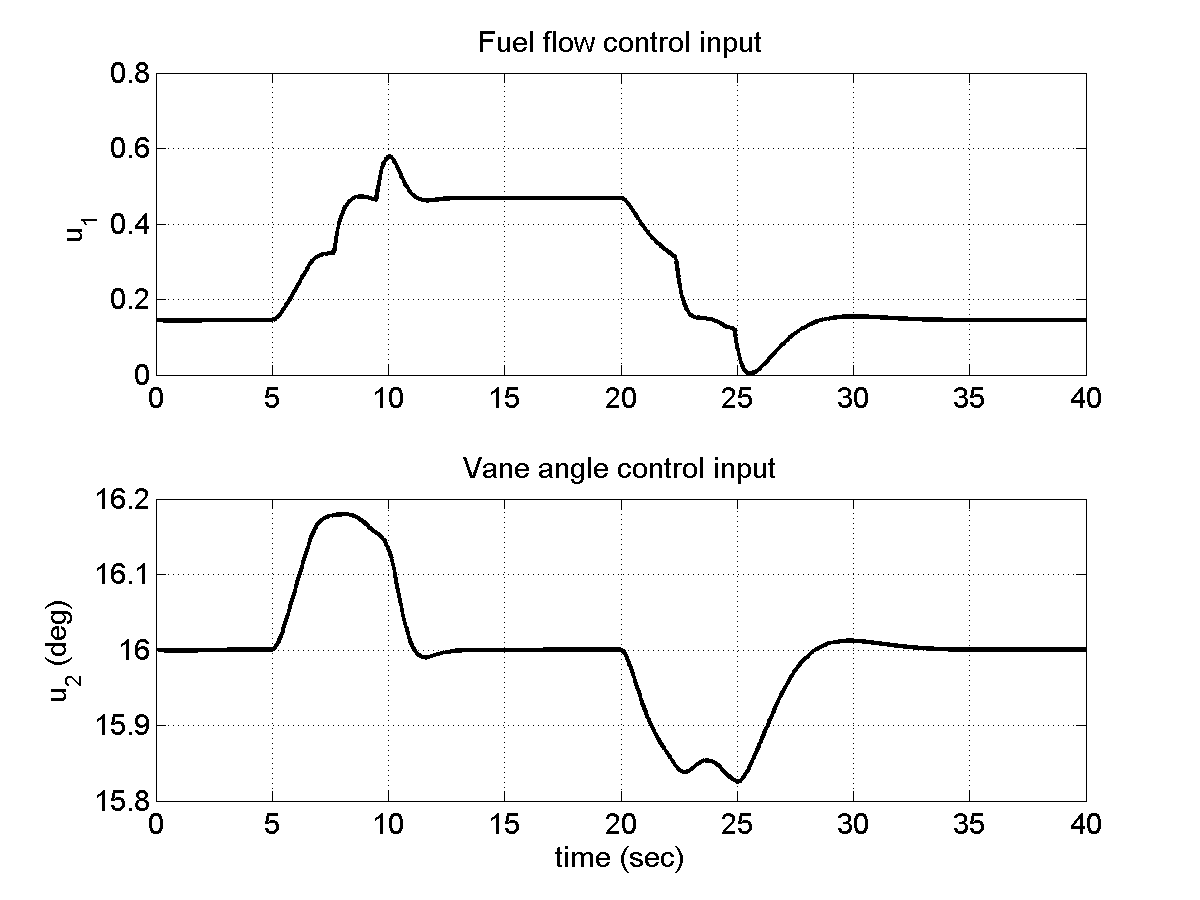}}
\caption{Fuel and prop pitch angle control inputs ($u(t)$)}\label{fig_cg_19}
\end{minipage}
\end{figure}

\begin{figure}[!ht]
\centering
\begin{minipage}[l]{3.2in}
\centering
\resizebox{3.2in}{!}{\includegraphics{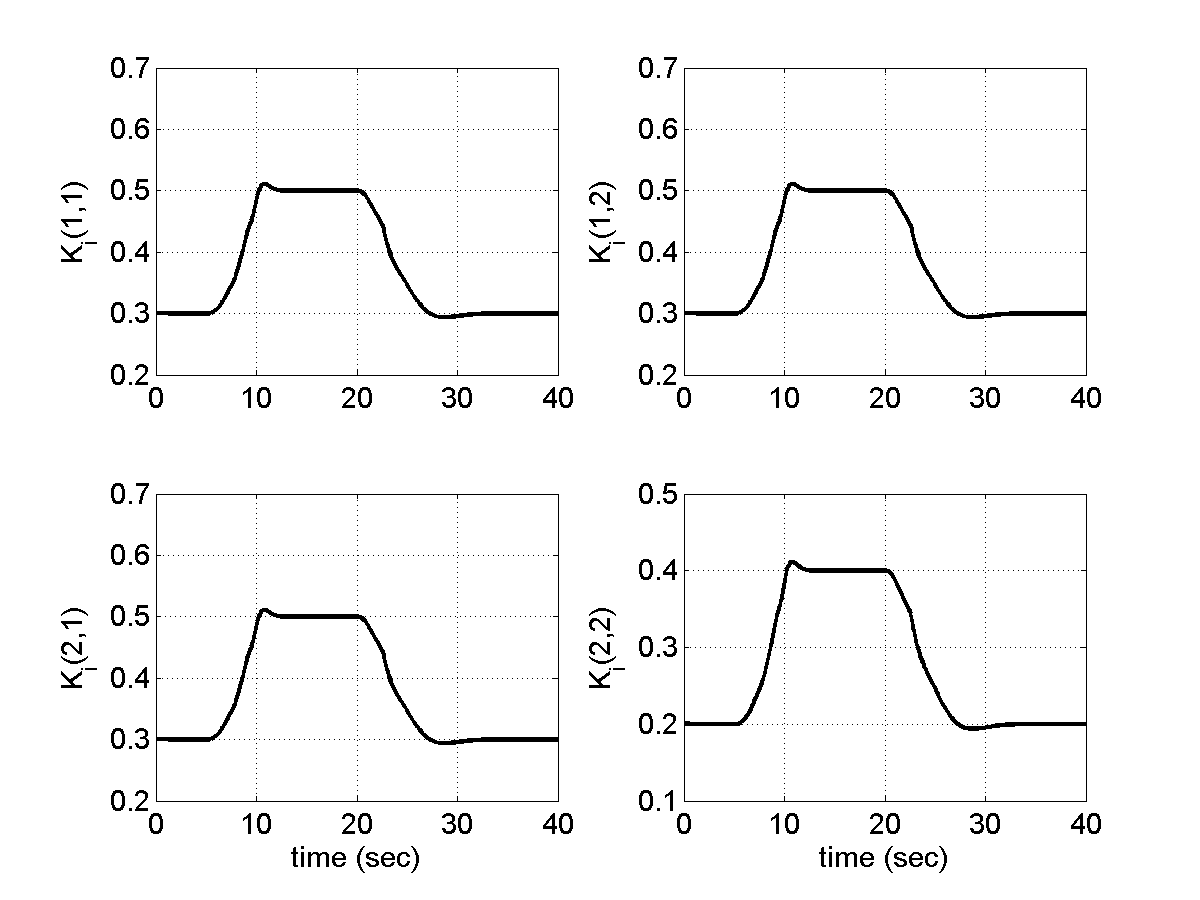}}
\caption{Controllers integral gain matrix ($K_i(\alpha)$) parameters histories}\label{fig_cg_20}
\end{minipage}
\begin{minipage}[r]{3.2in}
\centering
\resizebox{3.2in}{!}{\includegraphics{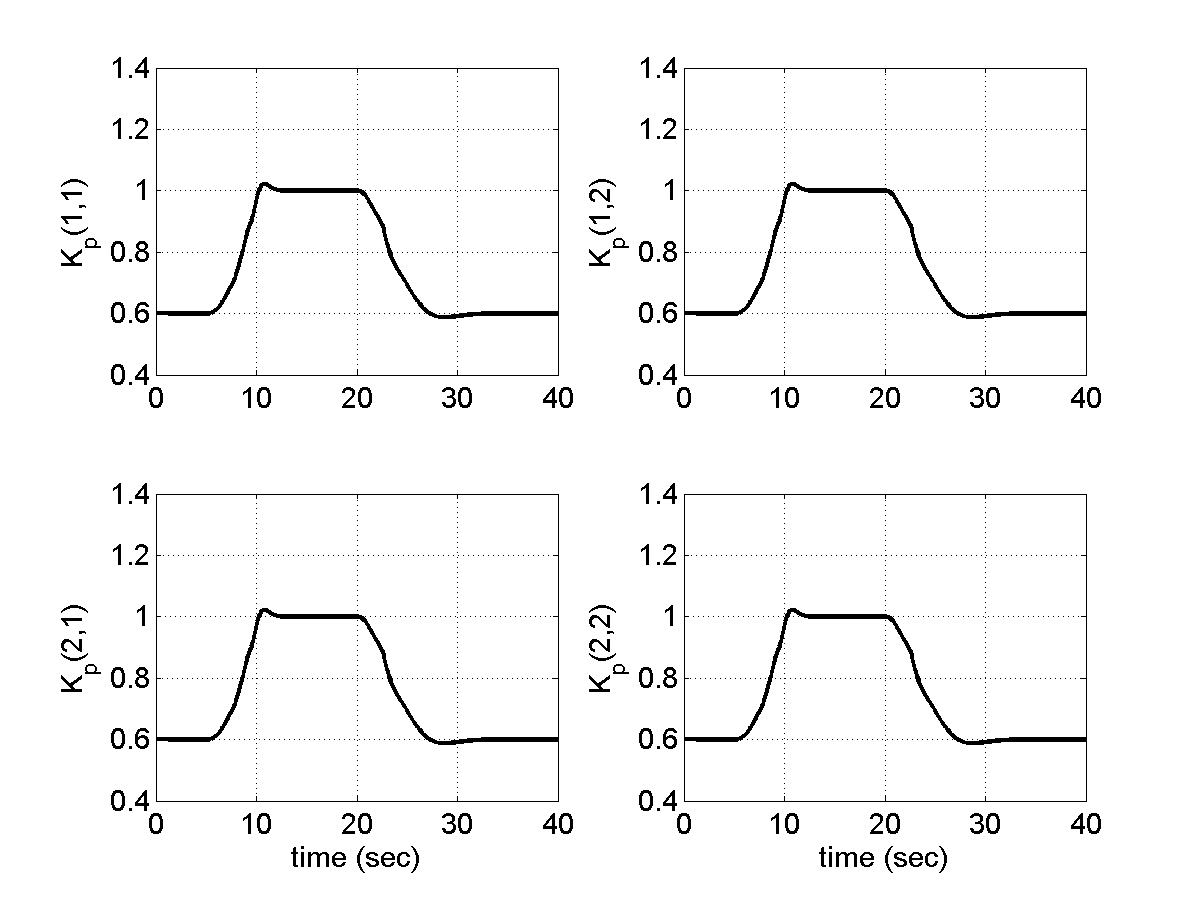}}
\caption{Controllers proportional gain matrix ($K_p(\alpha)$) parameters histories }\label{fig_cg_21}
\end{minipage}
\end{figure}

\begin{figure}[!ht]
\centering
\begin{minipage}[l]{3.2in}
\centering
\resizebox{3.2in}{!}{\includegraphics{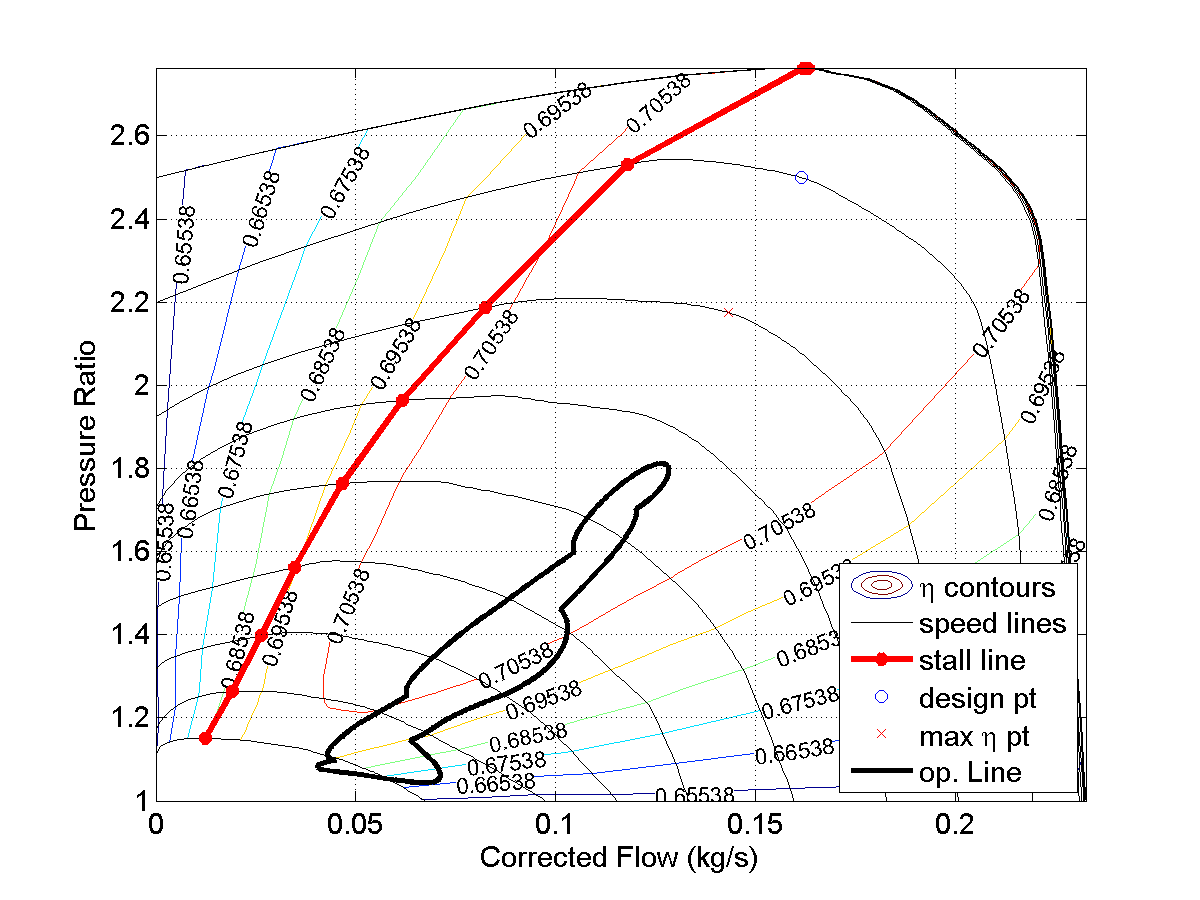}}
\caption{JetCat SPT5 engine compressor map with operating line}\label{fig_cg_22}
\end{minipage}
\begin{minipage}[r]{3.2in}
\centering
\resizebox{3.2in}{!}{\includegraphics{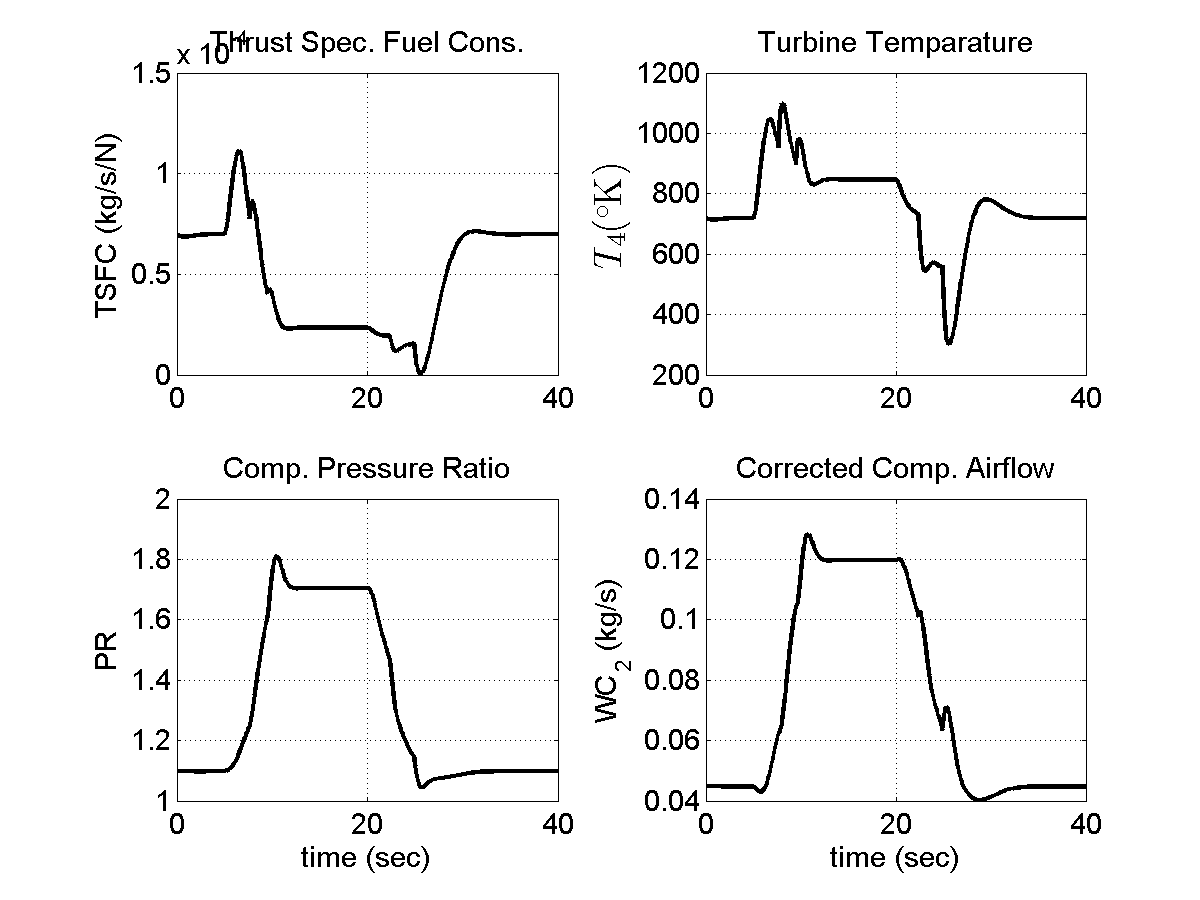}}
\caption{Turbine temperature, TSFC, compressor overall pressure ratio and air flow rate histories}\label{fig_cg_23}
\end{minipage}
\end{figure}

To show that the designed gain scheduled controller works properly on JetCat engine we used it to control the engine from idle to cruise condition and then again back to idle condition in a stable manner and with good performance. Simulation results are shown in figures \ref{fig_cg_01} to \ref{fig_cg_23}.

Figure \ref{fig_cg_01}, shows the history of the norm of closed-loop system matrix $||A_{cl}(t)||$, and its rate  $||\dot{A}_{cl}(t)||$. As it can be seen the figure shows the boundedness of these two variables in accordance with Assumption \ref{ass1} where $k_A=7.4539$, and $L_A=1.0106$. Figure \ref{fig_cg_02}, shows the history of the closed-loop system matrix eigenvalues $\lambda \{ A_{cl} \}$. As it is apparent, all the six eigenvalues remain negative with the time change of the scheduling parameter $\alpha$, and hence satisfies assumption \ref{ass2} of the stability theorem.

Figure \ref{fig_cg_03}, shows the history of the scheduling parameter which is $\alpha=p(y)=||y||=||x||$. It is also shows the history of the switching function, which is defined based on the norm of the spool speed equilibrium values vector. As it is apparent from the plot, engine operated in the vicinity of at least three equilibrium points to be able to accelerate from idle to cruise condition. The norm of the scheduling parameter rate $\dot{\alpha}(t)=\frac{x^T \dot{x}}{||x(t)||}$, also has been plotted. Figure \ref{fig_cg_04}, shows the history the norms of the output vector and its rate. This satisfies the condition of Theorem \ref{thm1} with $\epsilon_y=0.0025$. Using formulas from \cite{gainsched-shamma-1988}, we can compute $m=496.7476$, and $\lambda=0.5271$.

Figure \ref{fig_cg_08}, shows the history of the quadratic time varying Lyapunov function of the unforced closed loop system (\ref{eqn_gs50}). As it is apparent, $V(t)=\delta X^T P(t)\delta X$, is decrescent and bounded from above and below. The history of $\dot{V}$, shows that it is non-positive for all $t>0$, so the exponential stability of the slowly varying system (\ref{eqn_gs50}) with a gains-scheduling controller is guaranteed. Figure \ref{fig_cg_09}, shows the rate of change of the reference signals for the outputs of the system. The outputs in this simulation are core and spool speed. $||\dot{r}|| < 0.15$, which corresponds to the assumption of the $||\dot{r}||$ boundedness in the theorem \ref{thm3}.

Figure \ref{fig_cg_10}, shows the history of the plant states which are core and fan spool speeds. Figure \ref{fig_cg_11}, shows the time histories of the controller states. Figure \ref{fig_cg_12}, shows the phase plot for core and fan spool dynamics. Figure \ref{fig_cg_13}, shows the time history of the fan and core spool accelerations, i.e. $\dot{N}_1$ and $\dot{N}_2$.

 Figures \ref{fig_cg_14} and \ref{fig_cg_15}, show the core and fan spool speeds tracking their reference signals. Figure \ref{fig_cg_16}, shows the history of thrust and it is following its reference command from idle to cruise condition and then back to the idle for standard day, sea level condition. Figure \ref{fig_cg_17}, shows the control inputs to the augmented system, $v(t)=[v_1(t), v_2(t)]^T$, each element corresponding to one of the control inputs to the original system. To keep the engine dynamics within the limits of the operation for the available engine model, some limits have been defined on the augmented system control input, $-0.18 \leq v(t) \leq 0.18$. This limits will help to keep the fuel control input non-negative and also limits the rate of the control, $\dot{u}(t)$.

Figure \ref{fig_cg_18}, shows time rate of fuel and prop pitch angle inputs.  Figure \ref{fig_cg_19}, shows fuel flow and propeller pitch angle histories as control inputs. For better performance and also to keep the engine in the safe range of operation limits has been defined for the augmented control inputs. Figures \ref{fig_cg_20} and \ref{fig_cg_21}, show the baseline fuel controller integral ($K_i(\alpha)$) and proportional ($K_p(\alpha)$) gain matrices histories. These gains have been obtained by interpolation using the previously deigned fixed-gain controllers, each one corresponding to a equilibrium point of the engine. The numerical values of these gains are mentioned in equations (\ref{eqn_gs71}) to (\ref{eqn_gs73}).

Figure \ref{fig_cg_22}, shows JetCat SPT5 turboshaft engine compressor map. In this map the approximate stall line and also the operating line for this simulation has been shown. The engine operates in a safe region with a big stall margin during its acceleration from idle to cruise and again decelerating back to idle condition. Figure \ref{fig_cg_23}, shows the histories of turbine temperature, thrust specific fuel consumption (TSFC), compressor pressure ratio and corrected air flow rate.

Using gain-scheduling control technique, this simulation shows the possibility of controlling the linear parameter varying model of the turboshaft engine for a large transient throttle. This case study, simulates the engine accelerating from idle thrust to the cruise condition and then decelerates back to the idle condition in the standard day sea level condition. To tune the controller we need to trade off between the settling time and overshoot percentage. For large throttle transients, the existing controller regulates the controlled variable, additional limiters, will be developed to protect critical engine variables from exceeding physical bounds and to ensure safer operation.

%%%%%%%%%%%%%%%%%%%%%%%%%%%%%%%%%%%%%%%%%%%%%%%%%%%%%%%%%%%%%%%%%%%%%%
%\section{Conclusions and Future Works}
\section{Conclusions}

We developed a gain-scheduling controller with stability guarantees for nonlinear gas turbine engine systems. Using global linearization and LMI techniques, we showed guaranteed absolute stability for closed loop gas turbine engine systems with gain-scheduling controllers. Simulation results presented to show the applicability of the proposed controller to the nonlinear physics-based JetCat SPT5 turboshaft engine model for large transients from idle to cruise condition and vice versa.

%The future research directions are:
%\begin{itemize}
%	\item - Developing adaptive gain-scheduling controller for the engine, using the developed gain-scheduling controller to construct a time-varying reference model for the adaptive controller;
%	\item - Designing tools, such as a time-varying integrators \cite{gainsched-shamma-1988}, to relax the slowly-varying condition on the adaptive gain-scheduling controller;
%	\item - Developing a reference governor structure to handle the constraints on the system states and control inputs, based on the method presented in \cite{RefGovernor-bemporad-1998};
%	\item - Developing the decentralized version of the adaptive gain-scheduling controller to achieve plug-and-play control technology for gas turbine engines.
%\end{itemize}

%%%%%%%%%%%%%%%%%%%%%%%%%%%%%%%%%%%%%%%%%%%%%%%%%%%%%%%%%%%%%%%%%%%%%%
\section*{Acknowledgment}
This material is based upon the work supported by the Air Force Research Laboratory (AFRL).

%%%%%%%%%%%%%%%%%%%%%%%%%%%%%%%%%
\bibliographystyle{plain}
\bibliography{DistGainSchedulTurbineControl}

%%%%%%%%%%%%%%%%%%%%%%%%%%%%%%%%%%%%%%%%%%%%%%%%%%%%%%%%%%%%%%%%%%%%%%%
%\appendix       %%% starting appendix
%\section*{Appendix A: Head of First Appendix}
%Avoid Appendices if possible.

%%%%%%%%%%%%%%%%%%%%%%%%%%%%%%%%%%%%%%%%%%%%%%%%%%%%%%%%%%%%%%%%%%%%%%
%\end{multicols}

\end{document}